\documentclass[aps,pra,showpacs,preprint]{revtex4-1}

\usepackage{amsmath,amssymb,amsfonts,amsthm,}
\usepackage{graphicx}
\usepackage{graphicx}
\usepackage{txfonts}
\usepackage{wasysym}
\usepackage{stackrel}
\usepackage{stackengine}
\usepackage{hyperref}
\usepackage[T1]{fontenc}

\newcommand{\ieloperator}[1]{\; {#1}_{\scalebox{0.82}{\{}\!q\!\scalebox{0.82}{\}}} \;}
\newcommand{\oeloperator}[1]{\; {}_{\scalebox{0.82}{\{}\!q\!\scalebox{0.82}{\}} \!}{#1} \;}
\newcommand{\ileoperator}[1]{\; {#1}_{[q]} \;}
\newcommand{\oleoperator}[1]{\; {}_{[q]}{#1} \;}
\newcommand{\ioperator}[1]{\; {#1}_{\langle q \rangle} \;}
\newcommand{\ooperator}[1]{\; {}_{\langle q\rangle}{#1} \;}
\newcommand\owedge{\stackMath\mathbin{\stackinset{c}{0ex}{c}{0ex}{\scriptstyle{\wedge}}{\medcirc}}}
\newcommand\opm{\stackMath\mathbin{\stackinset{c}{0ex}{c}{0ex}{\scriptstyle{\pm}}{\bigcirc}}}

\begin{document}

\title{Deformed mathematical objects stemming from the $q$-logarithm function}

\author{Ernesto P. Borges} 
\email{ernesto@ufba.br} 
\affiliation{Instituto de Fisica, Universidade Federal da Bahia,
             Rua Barao de Jeremoabo, 40170-115 Salvador--BA, Brasil}
\author{Bruno G. da Costa}
\email{bruno.costa@ifsertao-pe.edu.br} 
\affiliation{Instituto Federal de Educa\c c\~ao, Ci\^encia e Tecnologia do
             Sert\~ao Pernambucano, {\it Campus} Petrolina, 
             BR 407, km 08, 56314-520 Petrolina, Pernambuco, Brazil}

\date{\today} 

\begin{abstract}
Generalized numbers, arithmetic operators and derivative operators,
grouped in four classes based on symmetry features,
are introduced.
Their building element is the pair of $q$-logarithm/$q$-exponential
inverse functions.
Some of the objects were previously described in the literature,
while others are newly defined.
Commutativity, associativity and distributivity,
and also a pair of linear/nonlinear derivatives
are observed within each class.
Two entropic functionals emerge from the formalism,
one of them is the nonadditive Tsallis entropy.
\end{abstract}

\maketitle 

\section{Introduction}

Extensivity of an entropy is expressed as $S$ being proportional 
to the number $N$ of elements of a composed system.
The hypervolume $\Omega$ of the phase space of a system composed by 
independent subsystems increases with 
the product of the hypervolumes $\mu_i$ of the corresponding 
subspaces of its elements ($\mu_i >1$).
For identical and independent subsystems,
the phase space exponentially increases with the number of elements,
$\Omega = \mu_1^N$,
and thus the Boltzmann entropy is proportional to $N$:
$S = k \ln \Omega = N k \ln \mu_1$,
i.e., it is extensive.
Correlations between subsystems make the hypervolume of the phase space 
smaller than that of the product of the hypervolumes of its subsystems,
and particular kinds of strong correlations make the phase space 
to asymptotically increase as a power law, 
a much slower rate than the exponential law;
in these cases the Boltzmann entropy is no longer extensive.
For such special cases, 
--- and there are plenty of observational, experimental and numerical examples ---,
the nonadditive entropy $S_q$ \cite{Tsallis-1988}
becomes proportional to $N$,
recovering extensivity, 
which is a central property for connecting with thermodynamics
(see details and further implications of extensivity 
 in Ref.\ \cite{Tsallis-springer-2009}).
The mathematical property that plays this role is 
a generalized multiplication operator defined in Ref.\ \cite{Borges-2004}.
The present paper identifies four classes of generalized algebras
associated with the nonextensive formalism in a broader point of view,
one of them contains the abovementioned generalized multiplication.
These developments hopefully help to understand the underlying mathematical
structures that support the nonextensive statistical mechanics.

The Tsallis nonadditive entropy $S_q$ 
has induced investigations on deformed mathematical structures aiming 
to represent relations of the nonextensive framework through expressions 
formally similar to the standard Boltzmann-Gibbs (BG) statistical mechanics.
The definition of the generalized logarithm function (the $q$-logarithm)
\cite{ct-quimicanova}
\begin{equation}
 \label{eq:q-logarithm}
 \displaystyle {\ln}_q\, x \equiv \frac{x^{1-q} -1}{1-q}, \quad (x > 0),
\end{equation}
allowed to rewrite $S_q \equiv k \, (q-1)^{-1} \, \big( 1-\sum_i^W p_i^q \big)$ 
(in its discrete version) as
\begin{eqnarray}
 \label{eq:Sq}
 \begin{array}{lll}
 \medskip
  S_q&=& - k \sum_i^W p_i^q \ln_q p_i,
\\
     &=&   k \sum_i^W p_i \ln_q (1/p_i)
 \end{array}
\end{eqnarray}
(sum over $W$ microstates, each one labeled $i$,
 with their corresponding probabilities $p_i$,
 $k$ is a positive constant, 
 $q \in \mathbb{R}$ is the generalizing entropic index).
Ordinary formalism is recovered as $q \to 1$ 
($\ln_1 x = \ln x$; $S_1 = S_{\text{BG}}$),
equiprobability yields $S_q[p_i=1/W]= k \ln_q W$.
The $q$-logarithm presents the limiting cases
\begin{eqnarray}
 \label{eq:qlog-limiting_cases}
 \lim_{x \to 0^+} \ln_q x & = & \left\{
                                     \begin{array}{ll}
                                      \frac{-1}{1-q}, & \quad q < 1,
                                      \\ \noalign{\smallskip}
                                      - \infty,       & \quad q \ge 1,
                                     \end{array}
                               \right.
 \\ \noalign{\smallskip}
 \lim_{x \to \infty} \ln_q x & = & \left\{
                                          \begin{array}{ll}
                                           \displaystyle
                                           \infty,        & \quad q \le 1,
                                           \\ \noalign{\smallskip}
                                           \frac{1}{q-1}, & \quad q > 1.
                                          \end{array}
                                   \right.
\end{eqnarray}
Its inverse, the $q$-exponential, is
\begin{equation}
 \label{eq:q-exponential}
  \exp_q(x) = 
        \left\{
               \begin{array}{ll}
                \left[1+(1-q)x\right]^\frac{1}{1-q} \,
                \theta\left(x+\frac{1}{1-q}\right),
                & q < 1,
                \\ \noalign{\smallskip}
                \text{e}^x, & q = 1, 
                \\ \noalign{\smallskip}
                \displaystyle
                \frac{1}{
                         \left[1-(q-1)x\right]^\frac{1}{q-1} \,
                         \theta\left( \frac{1}{q-1}-x \right)
                        }\:,
                & q > 1, 
              \end{array}
        \right.
\end{equation} 
($\theta(x)$ is the Heaviside step function)
that is more compactly written as
$\exp_q(x) = [1 + (1-q)x]_+^{1/(1-q)}$,
with the symbol $[\cdot]_+ \equiv \text{max}\{0,\cdot\}$, ---
the subscript symbol {\scriptsize $+$} encompasses the Heaviside function.
The Heaviside step function $\theta(x)$ defines the cutoff condition:
the $q$-exponential is set to zero for $q<1$ and $x<-1/(1-q)$,
and diverges for $q>1$ and $x>1/(q-1)$.
In the following we use either notations $\exp_q(x)$ or $\text{e}_q^x$, 
equivalently. 
Some properties of $q$-logarithm and $q$-exponential functions 
may be found in 
\cite{Yamano-2002,Naudts-2002,Tsallis-springer-2009,Naudts-springer-2011}.

The $q$-logarithm of a product splits into a nonadditive form for $q \ne 1$:
\begin{equation}
 \label{eq:nonadditivity-q-log}
 \ln_q (xy) = \ln_q x + \ln_q y + (1-q) \, \ln_q x \, \ln_q y.
\end{equation}
This property had triggered the definition of new generalized arithmetic 
operators:
(\textit{i}) 
what if the right hand side (r.h.s.) of this expression is viewed as
the definition of a generalized addition of $q$-logarithms?
Answer: Eq.\ (4) of \cite{Borges-2004}, 
Eq.\ (7) of \cite{lemans-2003},
Eq.\ (\ref{eq:ole-sum}) of the present work.
(\textit{ii}) 
What should be the argument of the $q$-logarithm of the 
left hand side (l.h.s.) of (\ref{eq:nonadditivity-q-log}) 
if its r.h.s.\ were an ordinary addition, 
instead of the generalized addition just defined?
Answer: Eq.\ (7) of \cite{Borges-2004}, 
Eq.\ (8) of \cite{lemans-2003},
Eq.\ (\ref{eq:oel-product}) of the present work.
Since then, these operators have been usually referred to as 
$q$-addition and $q$-multiplication, 
or, more colloquially, $q$-sum and $q$-product.
This $q$-multiplication is the one that makes $S_q$ extensive, 
mentioned previously, 
and it is not distributive with respect to the $q$-addition,
and Nivanen et al.\ \cite{lemans-2003} 
identified additional deformed operators, 
recovering the distributivity [their Eq.\ (24)--(28)].
In an extension of that work 
by the same authors with collaborators \cite{lemans-2009},
the $q$-multiplication and the $q$-addition were identified 
as belonging to two different classes, 
and further operators were defined.

Examples of mathematical developments along these lines:
spiral generalizations of the trigonometric and hyperbolic functions
through the extension of Euler's formula to the complex domain
\cite{Borges-1998},
generalization of derivative operators
\cite{Borges-2004,Nobre-RegoMonteiro-Tsallis-2011},
generalizations of Laplace and Fourier transforms,
representations of the Dirac delta function
\cite{Lenzi-Borges-Mendes-1999,
      Jauregui-Tsallis-2010,
      Sicuro-Tsallis-2017},
two parameter extensions for the logarithm and exponential
and their related algebras
\cite{Schwammle-Tsallis-2007,Cardoso-Borges-Lobao-Pinho-2008}
etc.
Other deformed mathematical structures had been introduced,
particularly the Kaniadakis formalism
\cite{Kaniadakis-physa-2001,
      Kaniadakis-pre-2002,
      Kaniadakis-Scarfone-2002,
      Kaniadakis-Lissia-Scarfone-2004},
from which some of the developments 
within the nonextensive context have been inspired.
Generalization of algebras has been recently proposed 
\cite{Gomez-Borges-g-algebra},
conforming to the group entropy theory
\cite{Tempesta-2011}.

Examples of physical systems described by nonextensive statistical mechanics:
anomalous diffusion of cold atoms in dissipative optical lattices
\cite{Lutz-Renzoni-2013},
anomalous diffusion in granular matter
\cite{Combe-Richefeu-Stasiak-Atman-2015},
experimental high energy physics
\cite{Wong-Wilk-2013},
observational high energy physics in cosmic rays
\cite{yalcin-beck-2018}.
An up-to-date bibliography may be found at the site \cite{temuco}.

The present paper revisits generalized algebras and calculus 
motivated by the nonextensive formalism in a broader point of view.
It identifies the basic arithmetic operators for four complementary classes, 
and defines a pair of linear/nonlinear derivative for each one.
A connection with entropic functionals is established.
The starting point is the definition of the generalized numbers.

The paper is organized as follows.
Section \ref{sec:q-numbers} introduces four deformed numbers,
by combining the pair of the inverse logarithm/exponential functions
and their generalized forms.
Section \ref{sec:q-algebras} explores each class of deformed arithmetics, 
derived from the generalized numbers.
Section \ref{sec:q-calculus} is dedicated to the deformed calculus
emerged from the infinitesimal deformed differences. 
Two possibilities are focused: a linear and a nonlinear deformed derivative.
A connection between these structures with entropic functionals
is addressed in section \ref{sec:entropy-generator}.
Particularly the nonadditive entropy $S_q$ is alternatively obtained 
through a procedure that uses one of the generalized powers 
defined in section \ref{sec:q-algebras}.
Section \ref{sec:conclusions} draws our final remarks 
and points towards new perspectives.
Throughout the text, many expressions use symbols 
designed for compactness.  
Some of them appear in their explicit forms in the Appendix.

\section{\label{sec:q-numbers}Deformed $q$-numbers}

Among the mathematical objects, one even more fundamental deserves attention
within this context, namely, the very concept of number.
This was implicitly advanced within the nonextensive formalism
in Ref.\ \cite{Borges-1998}, through the variable $\zeta_q = \ln \text{e}_q^z$ 
($z \in \mathbb{C}$)
used in the generalization of Euler's formula,
that may be read as a complex generalized number 
[see Eq.\ (22) of \cite{Borges-1998}].
Deformed numerical sets
($q$-natural $\mathbb{N}_q$, $q$-integer $\mathbb{Z}_q$, 
 $q$-rational $\mathbb{Q}_q$,
 $q$-real $\mathbb{R}_q$ numbers)
had been considered following Peano-like axioms
and generalized arithmetic operators had consistently been defined
\cite{Lobao-bjp-2009}.
Those generalized numbers are a transformation 
of the so-called {\it $Q$-analog of $n$} --- 
$Q$ standing for quantum, within the context of quantum calculus 
(we write it with upper case $Q$ to avoid confusion 
with the present lower case index $q$)
\cite{Kac-Cheung-2002}:
\begin{equation}
 \label{eq:quantum-number}
 [n]_Q = \frac{Q^n - 1}{Q-1},
\end{equation}
from which we borrow the idea of a $q$-number.
This connection had been previously realized, see \cite{Tsallis-1994}.
Deformation of reals had also been reported in Ref.\ \cite{Kalogeropoulos-2012}.

Given a continuous, analytical, monotonous, invertible function $f(x)$
generalized through a real parameter $q$ that recovers the ordinary case
as a limiting procedure (in this context, $q \to 1$),
we introduce the generalized numbers through four combinations 
such as the ordinary case is identically recovered:
\begin{subequations}
\label{eq:generalized-numbers}
\begin{eqnarray}
 \label{eq:iff-1}
      [x]_q \;\; &=& f        \big( f_q^{-1}(x) \big), \\ 
 \label{eq:off-1}
  {}_q[x] \;\;   &=& f_q      \big( f^{-1}  (x) \big), \\
 \label{eq:if-1f}
     \{x\}_q     &=& f^{-1}   \big( f_q     (x) \big), \\
 \label{eq:of-1f}
 {}_q\{x\} \;\;  &=& f_q^{-1} \big( f       (x) \big). 
\end{eqnarray}
\end{subequations}%

The adopted notation obeys the following criteria:
the square brackets are used when $f_q^{-1}$ (or $f^{-1}$) 
is the argument of $f$ (or $f_q$)
and the curly brackets are used when $f_q$ (or $f$) 
is the argument of $f^{-1}$ (or $f_q^{-1}$);
the function labeled as $f$ is arbitrary.
The deformation parameter $q$ is used as a subscripted {\it postfix}
if the {\it inner} function is deformed, referred to as i-number,
Eq.\ (\ref{eq:iff-1}) and (\ref{eq:if-1f}),
and as a subscripted {\it prefix}
if the {\it outer} function is deformed, referred to as o-number,
Eq.\ (\ref{eq:off-1}) and (\ref{eq:of-1f})
(in analogy with the notation employed for the 
 generalized hypergeometric series ---
 in that case, prefix for the numerator, postfix for the denominator).

The pair of i/o numbers are inverse of each other, and thus
\begin{equation}
 \label{eq:inverse}
   {}_q\big[\:      [x]_q\, \big] 
 =     \big[\,  {}_q[x]\:   \big]_q 
 = {}_q\big\{\:    \{x\}_q\,\big\} 
 =     \big\{\,{}_q\{x\}\:  \big\}_q 
 = x.
\end{equation}

To be more specific to the case we are focusing upon,
we define $f(x)=\ln x$, and, consequently, $f^{-1}(x)= \text{e}^x$.
It follows the le-numbers
(l stands for logarithm and e stands for exponential, 
 `le' expresses the order in which the functions are taken)
\begin{subequations}
  \label{eq:le-numbers}
\begin{eqnarray}
  \label{eq:ile-number}
   [x]_q  =& \ln \text{e}_q^x \, \quad (\text{ile-number}), 
  \\
  \label{eq:ole-number}
   {}_q[x] \;  =& \; \ln_q \text{e}^x \quad  (\text{ole-number}),
\end{eqnarray}
\end{subequations}%
and the el-numbers
\begin{subequations}
\label{eq:el-numbers-R+}
\begin{eqnarray}
  \label{eq:iel-number-complex}
   \{x\}_q     =& \; \text{e}^{\,\ln_q x} \quad (\text{iel-number}),
  \\
  \label{eq:oel-number-complex}
   {}_q\{x\} \;  =& \; \text{e}_q^{\ln x} \;\; \quad (\text{oel-number}).
\end{eqnarray}
\end{subequations}%
Equations (\ref{eq:el-numbers-R+})
are constrained to $x \in \mathbb{R}_+$.
This limitation can be overcome, allowing $x \in \mathbb{R}$,
in analogy to what was done in Ref.\ \cite{Tsallis-Queiros-2007},
by ad hoc redefining the el-numbers as
\begin{subequations}
 \label{eq:el-numbers}
 \begin{eqnarray}
  \label{eq:iel-number}
   \{x\}_q   =& \; \text{sign}(x) \, \text{e}^{\,\ln_q |x|} \quad \! 
   (\text{iel-number}),
  \\
  \label{eq:oel-number}
   {}_q\{x\} \; =& \; \text{sign}(x) \, \text{e}_q^{\ln |x|} \quad 
   (\text{oel-number}),
 \end{eqnarray}
\end{subequations}%
with $\text{sign}(x) = x / |x|$ and $\text{sign}(0) \equiv 0$.
The present work uses the el-numbers as defined by 
Eq.\ (\ref{eq:el-numbers}),
but expressions are easily rewritten in its simpler form 
(\ref{eq:el-numbers-R+})
by taking into consideration the restricted domain. 

\begin{figure}[htb]
 \begin{minipage}{0.47\linewidth}
  \includegraphics[width=0.8\linewidth,keepaspectratio,clip=]{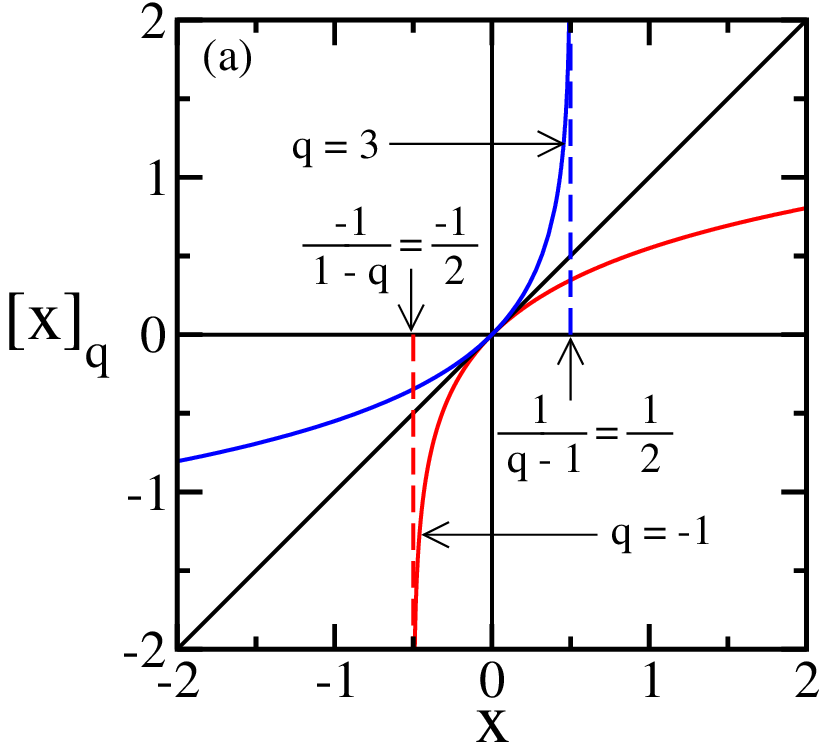}
 \end{minipage}
 \begin{minipage}{0.47\linewidth}
  \includegraphics[width=0.8\linewidth,keepaspectratio,clip=]{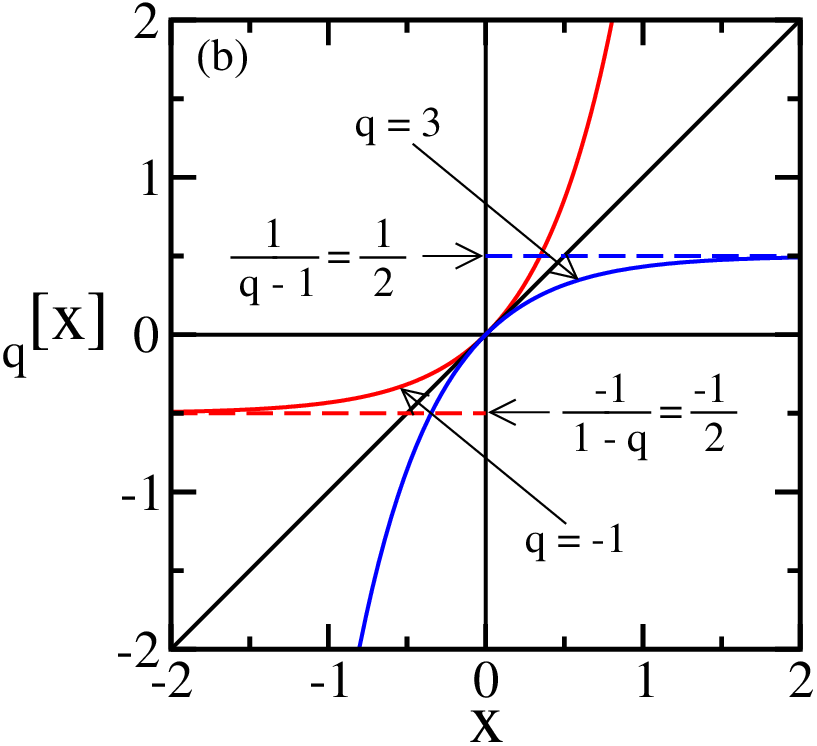}
 \end{minipage}
\\
 \begin{minipage}{0.47\linewidth}
  \includegraphics[width=0.8\linewidth,keepaspectratio,clip=]{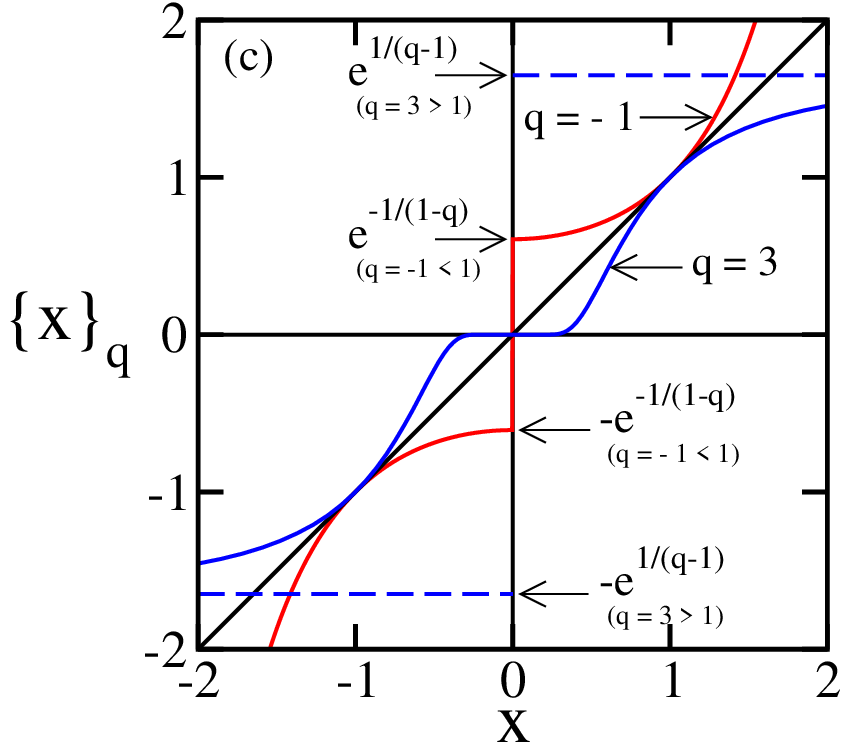}
 \end{minipage}
 \begin{minipage}{0.47\linewidth}
  \includegraphics[width=0.8\linewidth,keepaspectratio,clip=]{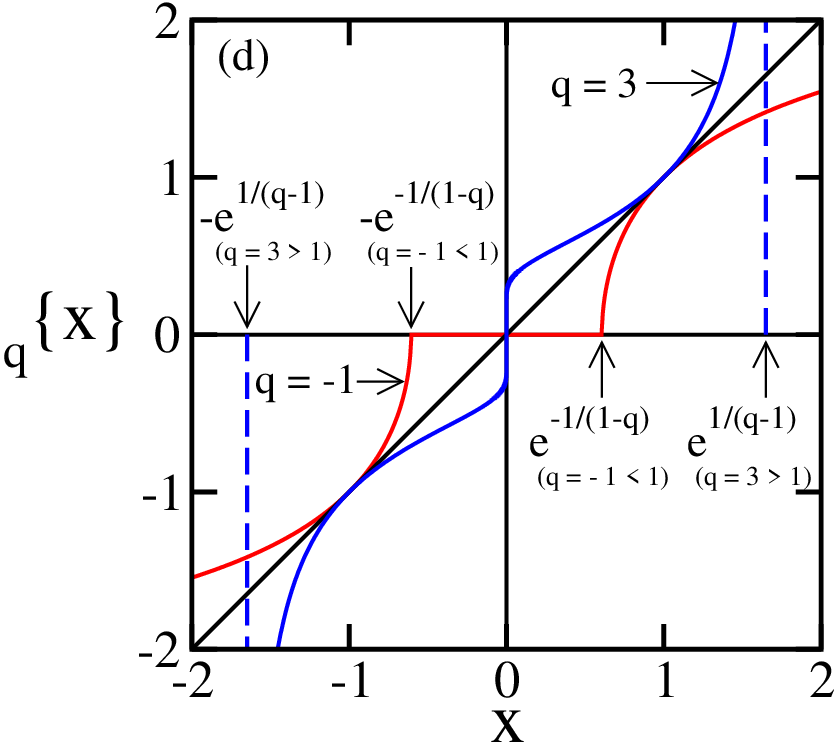}
 \end{minipage}
 \caption{$q$-numbers, illustrated with $q = -1$ (red), 1 (black), 3 (blue).
          (a) ile-number; 
              $[x \le -1/(1-q)]_{q<1} \to - \infty$, 
              illustrated by the vertical red asymptote for $q=-1$;
              $[x \ge 1/(q-1)]_{q>1} \to  \infty$,
              illustrated by the vertical blue asymptote for $q=3$.
          (b) ole-number;
              $\lim_{x \to -\infty} {}_{q<1}[x] = -1/(1-q)$,
              illustrated by the horizontal red asymptote for $q=-1$:
              $\lim_{x \to  \infty} {}_{q>1}[x] = 1/(q-1)$,
              illustrated by the horizontal blue asymptote for $q=3$.
          (c) iel-number;
              $\lim_{x \to 0^{\pm}} \{x\}_{q<1} = \pm \text{e}^{-1/(1-q)}$;
              illustrated for $q=-1$;
              $\lim_{x \to \pm \infty} \{x\}_{q>1} = \pm \text{e}^{1/(q-1)}$,
              illustrated by the horizontal blue asymptotes for $q=3$;
          (d) oel-number;
              ${}_{q<1}\{|x| \le \text{e}^{-1/(1-q)}\} = 0$,
              illustrated for $q=-1$;
              ${}_{q>1}\{|x| \ge \text{e}^{1/(q-1)}\} \to \text{sign}(x)\infty$,
              illustrated by the vertical blue asymptotes for $q=3$.
          \label{fig:q-numbers}
         }
\end{figure}

The le-numbers have one fixed point ($[x]=x$)
at $[0]_q=0$, and ${}_q[0]=0$
(ile and ole, respectively) for all values of $q \ne 1$.
The iel-numbers have two fixed points ($\{x\}=x$)
for $q<1$, at $\{\pm 1\}_q = \pm 1$, 
--- zero is not a fixed point for iel-numbers, since $\nexists \, \{0\}_q$
    ($\lim_{x \to 0^-} \{x\}_q = -\text{e}^{-1/(1-q)}$,
     $\lim_{x \to 0^+} \{x\}_q =  \text{e}^{-1/(1-q)}$) ---,
and three fixed points for $q \ge 1$, 
at $\{0\}_q = 0$ and $\{\pm 1\}_q = \pm 1$.
The oel-numbers have three fixed points, 
at ${}_q\{0\}=0$, and ${}_q\{\pm 1\}= \pm 1$. 
Due to the cutoff condition of the $q$-exponential,
${}_{q < 1}\big\{|x|<\text{e}^{1/(q-1)}\big\}=0$,
and due to the absolute values, the el-numbers are odd,
for both i and o deformed numbers ($\{-x\}=-\{x\}$).
le-numbers and el-numbers are monotonous crescent with the ordinary numbers,
i.e., if $x>y$, $[x] > [y]$ and $\{x\} > \{y\}$ 
for both i and o deformed numbers.
Exception may apply for oel-numbers:
it may happen $x>y$ but ${}_q\{x\}={}_{q}\{y\}=0$ for $q<1$
within the cutoff region, 
$|x| \le \exp\big(\!-1/(1-q)\big)$ and $|y| \le \exp\big(\!-1/(1-q)\big)$.
The inverse relations between ile/ole and iel/oel numbers 
expressed by Eq.\ (\ref{eq:inverse}) are valid outside the cutoff regions.
Figure \ref{fig:q-numbers} illustrates the four $q$-numbers.
These deformed numbers also satisfy the identities
\begin{subequations}
 \label{eq:qlog-qnumbers}
  \begin{eqnarray}
       \label{eq:ln-oelx}
       \big[\!\ln x   \big]_q  &=& \ln\big(   {}_q\{x\}   \,\big),     
\\
   {}_q\big[\!\ln x   \big]    &=& \ln\big(\,       \{x\}_q \,\big) \;\;\: 
                                = \;\; \ln_q x,
   \label{eq:qlog-as-log}
\\
       \big[\!\ln_q x \big]_q &=& \ln_q\big( {}_q\{x\}   \,\big) \;\; 
                               = \;\; \ln x,
\\
   {}_q\big[\!\ln_q x \big]   &=& \ln_q\big(\,     \{x\}_q \,\big),    
  \end{eqnarray}
\end{subequations}%
\begin{subequations}
 \label{eq:qexp-qnumbers}
  \begin{eqnarray}
       \label{eq:exp-olex}
       \big\{ \exp x   \,\big\}_q &=& \exp\big(\, {}_q[x] \,\big), 
\\
       \label{eq:exp-ilex}
   {}_q\big\{ \exp x   \,\big\}   &=& \exp\big(\;\:   [x]_q \,\big) \; 
                                   = \;\, \exp_q x,
\\
       \label{eq:expq-ilex}
       \big\{ \exp_q x \,\big\}_q &=& \exp_q\big(\, {}_q[x] \,\big)   \;
                                   =  \;\, \exp x,
\\
       \label{eq:expq-olex}
   {}_q\big\{ \exp_q x \,\big\}   &=& \exp_q\big(\:   [x]_q \,\big).
  \end{eqnarray}
\end{subequations}%

Whenever convenient and not ambiguous,
for the sake of compactness of notation, 
we henceforth may occasionally use the symbols
$\langle x \rangle_q$ to denote the i-numbers 
(either $[x]_q$ or $\{x\}_q$),
and
${}_q\langle x \rangle$ to denote the o-numbers 
(either ${}_q[x]$ or ${}_q\{x\}$),
and the most general case $\langle x \rangle$, without subscripts,
to denote any of the four generalized numbers
(not to be confound with mean value or the bra-ket symbols).
The expressions `generalized number' and `generalized variable'
are used interchangeably, just as the convenience of the context,
without restricting ourselves to rigorous mathematical distinction 
these concepts may have.

Next sections explore connections of these deformed numbers 
with their corresponding arithmetics and calculus.

\section{\label{sec:q-algebras}Deformed $q$-arithmetics}

Starting from the generalized numbers 
(\ref{eq:le-numbers}) and (\ref{eq:el-numbers})
we identify four generalized classes of arithmetics.
In this paper, the designation $q$-addition, $q$-multiplication etc.\ 
are ambiguous, and thus we introduce a different notation: 
the ile-, ole-, iel-, oel- arithmetic operators.
Particularly, and partially anticipating the results of the next subsections,
the deformed addition and subtraction of Ref.\ \cite{Borges-2004}
belong to the ole-algebra
(here symbolized by $\!\oleoperator{\oplus}$ and $\!\oleoperator{\ominus}$),
considered in subsection \ref{subsec:ole-algebra},
and the deformed multiplication and division of Ref.\ \cite{Borges-2004}
belong to the oel-algebra
(here symbolized by $\!\oeloperator{\otimes}$ and $\!\oeloperator{\oslash}$),
considered in subsection \ref{subsec:oel-algebra}.
By $q$-arithmetics we generically denote the set of the four arithmetics
described in this paper.
They can also be referred to as $q$-algebras, 
understood as algebras over the real numbers, or some subset of the reals.

An i-arithmetic operator is defined as the i-number of the ordinary
arithmetic operator of the corresponding o-numbers,
and, complementary,
an o-arithmetic operator is defined as the o-number of the ordinary
arithmetic operator of the corresponding i-numbers.
The generating rules follow the lines 
of the $\kappa$-arithmetic operators of Kaniadakis
\cite{Kaniadakis-physa-2001,Kaniadakis-pre-2002,Kaniadakis-Scarfone-2002},
more generally expressed by Eq.\ (1) of \cite{Jizba-2020} 
(also in \cite{Gomez-Borges-g-algebra}),
and are
\begin{subequations}
\label{eq:generating-rules}
 \begin{eqnarray}
 \label{eq:generating-rule-i}
  \text{{i-arithmetics:}} \qquad &
   x \ioperator{\Circle} y = \;\; \left\langle \: 
                                    {}_q\langle x \rangle \; 
                                   \circ \; 
                                    {}_q\langle y \rangle \: 
                                  \right\rangle_q,
 \end{eqnarray}
 \begin{eqnarray}
 \label{eq:generating-rule-o}
\text{{o-arithmetics:}} \qquad &
   x \ooperator{\Circle} y = \;   {}_q\left\langle \: 
                                      \langle x\rangle_q \; 
                                   \circ \; 
                                      \langle y \rangle_q \: 
                                      \right\rangle.
 \end{eqnarray}
\end{subequations}
 The symbol $\circ$, a small circle without subscripts,
 represents any general usual arithmetic operator,
 $\circ \in \{+,-,\times,\slash\}$;
 its generalized version is represented by a larger circle \Circle, 
 with bracket subscripts: 
 prefixed/postfixed, square/curly, in consonance with the case.
To avoid ambiguity in notation, 
the generalized {\it operators} are represented within a circle 
with their subscripts within brackets.
The generalized {\it numbers} are represented within brackets, 
with their subscripts without brackets%

Some general relations are valid for all cases
(the symbol $N$ without subscript generically represents the neutral
element of the addition for any of the four arithmetics
$N \in \{N_{\scalebox{0.57}{[+]}},\,
{}_{\scalebox{0.57}{[+]}}N,\,
N_{\scalebox{.57}{\{+\}}},\,
{}_{\scalebox{.57}{\{+\}}}\!N\}$;
similarly to $I$, the neutral element of the multiplication;
$A$, the absorbing element of the multiplication):
the neutral element of the deformed addition $N$,
such as
$x \, \oplus \, N = x$,
is the the corresponding deformed zero
(%
$N_{\scalebox{0.57}{[+]}} = [0]_q$,
${}_{\scalebox{0.57}{[+]}}N = {}_q[0]$,
$ N_{\scalebox{.57}{\{+\}}} = \{0\}_q $, 
$ {}_{\scalebox{.57}{\{+\}}}\!N = {}_q\{0\}$%
);
the deformed additive opposite of $x$, written as 
$\ominus \: x \equiv 0\ominus x$,
such that
$x \oplus (\ominus x) = N$,
and 
$x \ominus y = x \oplus (\ominus \: y)$.
Similarly, the neutral element of the deformed multiplication $I$,
$x \otimes I = x$,
is the corresponding deformed unity
($I_{\scalebox{.57}{[$\times$]}} = [1]_q$,
 ${}_{\scalebox{0.57}{[$\times$]}}I = {}_q[1]$,
 $I_{\scalebox{0.57}{\{$\times$\}}} = \{1\}_q$, 
 ${}_{\scalebox{0.57}{\{$\times$\}}}\!I = {}_q\{1\}$).
The deformed multiplicative inverse of $x$, written as $I \oslash x$,
is such that
$x \otimes (I \oslash x) = I$,
and
$x \oslash y = x \otimes (I \oslash y)$.
The absorbing element of the deformed multiplication $A$,
such that
$x \otimes A = A$,
coincides with the neutral element of the corresponding deformed addition $N$
($A_{\scalebox{0.57}{[$\times$]}} = N_{\scalebox{0.57}{[+]}} $,
 ${}_{\scalebox{0.57}{[$\times$]}}A = {}_{\scalebox{0.57}{[+]}}N $,
 $ A_{\scalebox{.57}{\{$\times$\}}} = N_{\scalebox{.57}{\{+\}}}$,
 $ {}_{\scalebox{.57}{\{$\times$\}}}\!A = {}_{\scalebox{.57}{\{+\}}}\!N$).
The deformed addition and multiplication are commutative
($x \oplus y = y \oplus x$,
 \ 
 $x \otimes y = y \otimes x$),
associative
[$x \oplus (y \oplus z) = (x \oplus y) \oplus z$,
 \
 $x \otimes (y \otimes z) = (x \otimes y) \otimes z$],
and the deformed multiplication is 
(left and right) distributive with respect to the deformed addition
\big[$x \otimes (y \oplus z) = (x \otimes y) \oplus (x \otimes z)$,
 \ 
 $(y \oplus z) \otimes x = (y \otimes x) \oplus (z \otimes x)$\big]
\cite{Jizba-2020}.
Some constraints may apply to these relations 
according to the case, to be detailed in the next subsections.

\subsection{\label{subsec:ile-algebra}ile-Arithmetics}

The ile-algebraic operators follow from the generating rule 
expressed by (\ref{eq:generating-rule-i}).
The ile-addition is
\begin{eqnarray}
 \label{eq:ile-sum}
 \begin{array}{lll}
 \smallskip
 x \ileoperator{\oplus} y &=& \left[ \, {}_q[x] \; + \; {}_q[y] \: \right]_q,
 \\
                        &=& \ln \exp_q \left( \, 
                                         \ln_q \text{e}^x + \ln_q \text{e}^y \, 
                                       \right).
 \end{array}
\end{eqnarray}
The neutral element of the ile-addition is 
$N_{\scalebox{0.57}{[+]}} = [0]_q = 0$,
and consequently the opposite ile-additive of $y$ is
\begin{eqnarray}
 \label{eq:ile-opposite}
    \ileoperator{\ominus} y = 
    \frac{1}{1-q} \ln \left( 2-\text{e}^{(1-q)y} \right).
\end{eqnarray}
The ile-difference (\ref{eq:generating-rule-i}) with 
$\ioperator{\Circle}=\ileoperator{\ominus}$,
\begin{eqnarray}
 \label{eq:ile-difference}
 \begin{array}{lll}
 \smallskip
 x \ileoperator{\ominus} y &=& \;\; \left[ \, 
                                           {}_q[x] \; - \; {}_q[y] \: 
                                    \right]_q,
 \\
                        &=& \ln \exp_q \left( \, 
                                         \ln_q \text{e}^x - \ln_q \text{e}^y \, 
                                       \right),
 \end{array}
\end{eqnarray}
consistently satisfies
$x \ileoperator{\ominus} y = 
   x \ileoperator{\oplus} ( \ileoperator{\ominus} y \:)$
for all $q$.
The ile-multiplication is
\begin{eqnarray}
 \label{eq:ile-product}
 \begin{array}{lll}
 \smallskip
 x \ileoperator{\otimes} y &=& \;\; \left[ \, 
                                          {}_q[x] \; \; {}_q[y] \: 
                                    \right]_q,
 \\
                        &=& \ln \exp_q \left( \, 
                                        \ln_q \text{e}^x \; \ln_q \text{e}^y \, 
                                       \right),
 \end{array}
\end{eqnarray}
with its neutral ile-multiplicative element 
$I_{\scalebox{.57}{[$\times$]}} = [1]_q = (1-q)^{-1}\ln(2-q)$ 
for $q<2$
($[1]_{q \ne 1} \ne 1$,
 $\nexists  I_{\scalebox{.57}{[$\times$]}}$ for $q \ge 2$),
and its ile-absorbing element
$A_{\scalebox{0.57}{[$\times$]}} = [0]_q = 0$,
for all $q$.
The ile-division is
\begin{eqnarray}
 \label{eq:ile-ratio}
 \begin{array}{lll}
 \smallskip
 x \ileoperator{\oslash} y &=& \;\; \left[\,
                                          \displaystyle
                                          \frac{{}_q[x]}{ {}_q[y]}\:
                                    \right]_q,
 \\
                        &=& \ln \exp_q \left(\displaystyle
                                             \frac{\ln_q \text{e}^x}
                                                  {\ln_q \text{e}^y}
                                       \right),
 \end{array}
\end{eqnarray}
and
$\nexists \; I_{\scalebox{0.57}{[$\times$]}} \ileoperator{\oslash} 0$.

The ile-power of $x$ is defined as the ile-multiplication of 
$n$ identical factors $x$,
\begin{eqnarray}
 \label{eq:ile-repeated-product}
  x \ileoperator{\owedge} n =
                                \ileoperator{\displaystyle\prod^n} x
                              = \left[ \;
                                       \left( 
                                             {}_q[x] \;
                                       \right)^n \,
                                \right]_q.
\end{eqnarray}
Its analytical extension from $n \in \mathbb{N}$ to $y \in \mathbb{R}$ 
is written as
\begin{equation}
 \label{eq:ile-power}
  x {\ileoperator{\owedge} y} = \ln \exp_q \left(\, 
                                                 (\ln_q \text{e}^x)^y \, 
                                           \right),
  \quad (x>0),
\end{equation}
with the particular cases:
$x{\ileoperator{\owedge} 0} = [1]_q$ ($x \ne 0$), \,
$x{\ileoperator{\owedge} 1} = x$ ($x \ne 0$), \,
$1\!{\ileoperator{\owedge} y} \ne 1$ (for $q \ne 1$), \,
$\lim_{x \to 0^+} (x{\ileoperator{\owedge} y}) = 0$ ($y > 0$), \,
$\lim_{x \to 0^+} (x{\ileoperator{\owedge} y}) \to \infty$ ($y < 0$), \,
and the trivial case
$x\,{\owedge_{[1]}\; y} = x^y$.
The ile-power is right-distributive with respect to the ile-multiplication:
$
 (x \ileoperator{\otimes} y) \ileoperator{\owedge} z 
 =
 (x \ileoperator{\owedge} z)
 \ileoperator{\otimes}
 (y \ileoperator{\owedge} z).
$

The repeated generalized addition defines a different 
generalized multiplication that can be named as dot-multiplication, 
identified by the symbol $\odot$,
to distinguish it from the previous generalized multiplication
(or times-multiplication), symbolized by $\otimes$
[Eq.\ (\ref{eq:ile-product}) for the ile class].
The repeated ile-addition is given by
\begin{eqnarray}
 \label{eq:ile-repeated-sum}
  \begin{array}{lll}
  \smallskip
  n \ileoperator{\odot} y  &=& \ileoperator{\displaystyle\sum^n} y,
  \\ \smallskip
                           &=& \bigg[\;
                                     \displaystyle\sum^n \; {}_q[y]\;
                               \bigg]_q,
  \\
                           &=& \ln \exp_q \left( n\:\ln_q \text{e}^y \right),
  \end{array}
\end{eqnarray}
where we have used the generalized summation symbol for the ile class,
$\ileoperator{\sum}$,
compatible with the notation adopted in this work.
Analytical extension from $n \in \mathbb{N}$ to $x \in \mathbb{R}$ yields
the non commutative generalized ile-dot-multiplication:
\begin{eqnarray}
 \label{eq:ile-dot-product}
  x \ileoperator{\odot} y  &=& 
        \frac{1}{1-q} \ln \left( x\, \text{e}^{(1-q)y} - (x-1) \right)_{\!+}.
\end{eqnarray}
The dot-multiplication with the unity has two behaviors, 
due to its non-commutativity.
The trivial case ($1 \odot y = y$) holds for the four classes
(for the ile-dot-multiplication of this subsection,
 as well as for the ole-, iel-, oel- of the subsections to come).
The other case, $x \odot 1$, connects the dot-multiplication 
with the deformed numbers. The ile-dot-multiplication with unity results
$x \ileoperator{\odot} 1 = \big[\, x \; {}_q[1] \,\big]_q$,
with ${}_q[1] = \ln_q e = \big(\text{e}^{1-q}-1\big)/(1-q) \ne 1$ for $q \ne 1$.
Repeated ile-dot-multiplication defines ile-dot-power,
not explicitly shown here.

\subsection{\label{subsec:ole-algebra}ole-Arithmetics}

The generating rule (\ref{eq:generating-rule-o}) 
defines the ole-algebraic operators.
The ole-addition (or ole-sum) is
\begin{eqnarray}
 \label{eq:ole-sum}
 \begin{array}{lll}
  \smallskip
  x \oleoperator{\oplus} y &=&  {}_{q}\left[ \, [x]_q \; + \; [y]_q \, \right],
  \\ \smallskip
                           &=& \ln_q \exp \left(\,
                                            \ln \text{e}_q^x + \ln \text{e}_q^y
                                          \right),
  \\
                           &=& x + y + (1-q)xy.
 \end{array}
\end{eqnarray}
Its neutral ole-additive element is ${}_{\scalebox{0.57}{[+]}}N = {}_q[0] = 0$
and the opposite ole-additive element $\oleoperator{\ominus} y$ such as
$y \oleoperator{\oplus} (0 \oleoperator{\ominus} y) = 0$ is
\begin{eqnarray}
 \label{eq:ole-opposite}
 \begin{array}{lll}
  \oleoperator{\ominus} y 
                          &=& \displaystyle\frac{-y}{1+(1-q)y},
                              \quad 
                              \big(y \ne 1/(q-1)\big),
 \end{array}
\end{eqnarray}
and, consequently, the ole-subtraction is
\begin{eqnarray}
 \label{eq:ole-difference}
 \begin{array}{lll}
  \smallskip
  x \oleoperator{\ominus} y &=&  {}_{q}\left[ \, 
                                             [x]_q \; - \; [y]_q \, 
                                       \right],
  \\ \smallskip
                           &=& \displaystyle
                               \ln_q \exp\left( \, 
                                          \ln \text{e}_q^x - \ln \text{e}_q^y \,
                                         \right),
  \\ 
                           &=& \displaystyle\frac{x-y}{1+(1-q)y}
 \end{array}
\end{eqnarray}
provided $y \ne 1/(q-1)$.
These are the generalized addition and subtraction of Ref.\ \cite{Borges-2004},
referred to as $q$-sum and $q$-difference, respectively
(see also Subsec.\ 3.3.3 of Ref.\ \cite{Tsallis-springer-2009}).

From (\ref{eq:generating-rule-o}),
the ole-product
\begin{eqnarray}
 \label{eq:ole-product}
 \begin{array}{lll}
 \smallskip
  x \oleoperator{\otimes} y &=&  {}_{q}\left[ \, 
                                              [x]_q \; [y]_q \: 
                                       \right],
  \\
                           &=& \displaystyle
                               \ln_q \exp\left(
                                           \ln \text{e}_q^x \; \ln \text{e}_q^y 
                                         \right),
 \end{array}
\end{eqnarray}
and its neutral ole-multiplicative element
${}_{\scalebox{0.57}{[$\times$]}}I = {}_q[1] 
  = \displaystyle\frac{\text{e}^{1-q}-1}{1-q} \ne 1$
for $q \ne 1$, together with the ole-division, 
\begin{eqnarray}
 \label{eq:ole-ratio}
 \begin{array}{lll}
  \smallskip
  x \oleoperator{\oslash} y &=&
             \stackrel[ {\stackrel[ \scalebox{0.8}{\textit{q}} ]{}{}} ]{}{} \!\!
                                \left[
                                      \displaystyle\frac{[x]_q}{[y]_q}
                                \right],
  \\
                            &=& \displaystyle
                                 \ln_q \exp \left(
                                                  \frac{\ln \text{e}_q^x}
                                                       {\ln \text{e}_q^y}
                                            \right),
 \end{array}
\end{eqnarray}
are coherent with the ole-multiplicative inverse element
${}_{\scalebox{0.57}{[$\times$]}}I\!\oleoperator{\oslash}y =
    \ln_q\exp\left(\,(\ln \text{e}_q^y)^{-1}\right)$.
The ole-absorbing element is 
${}_{\scalebox{0.57}{[$\times$]}}A = {}_q[0]=0$.
The generalized diamond multiplication defined 
by Eq.\ (24) of Ref.\ \cite{Lobao-bjp-2009}
is related to the ole-multiplication as
$x \oleoperator{\otimes} y = 
   (x \, \Diamond_{\!q} \, y) \oleoperator{\otimes} 1$,
and this expression connects the distributivity property
of the diamond multiplication with respect to the ole-addition
(Eq.\ (28) of Ref.\ \cite{Lobao-bjp-2009})
and the distributivity of the ole-multiplication with respect to 
this generalized addition.

The ole-power (the repeated ole-multiplication),
\begin{eqnarray}
 \label{eq:ole-repeated-product}
  x {\oleoperator{\owedge} n} =
                                \stackrel[ \scalebox{0.7}{\{\textit{q}\}} ]
                                         {}{} \!\!
                                \displaystyle\prod^n x
                              = {}_q\left[ \:
                                           \left(
                                                 [x]_q
                                           \right)^n
                                    \right],
\end{eqnarray}
after analytic continuation, becomes
\begin{equation}
 \label{eq:ole-power}
  x {\oleoperator{\owedge} y} = \ln_q \exp\left(\ln \text{e}_q^x\right)^y,
  \quad (x>0),
\end{equation}
with
$x {\!\oleoperator{\owedge} 0} = {}_q[1]$ ($x \ne 0$), \,
$x {\!\oleoperator{\owedge} 1} = x$ ($x \ne 0$), \,
$1 {\!\oleoperator{\owedge} y} \ne 1$ (for $q \ne 1$), \,
$\lim_{x\to 0^+} (x {\!\oleoperator{\owedge} y}) = 0$ ($y > 0$), \,
$\lim_{x\to 0^+} (x {\!\oleoperator{\owedge} y}) \to \infty$ 
 ($y < 0$ and $q<1$), \,
$\lim_{x\to 0^+} (x {\!\oleoperator{\owedge} y}) =1/(q-1)$ 
 ($y < 0$ and $q>1$), \,
$x\:{{}_{[1]}{\owedge} \; y} = x^y$.
The ole-power is right-distributive with respect to the ole-multiplication:
$
 (x \oleoperator{\otimes} y) \oleoperator{\owedge} z 
 =
 (x \oleoperator{\owedge} z)
 \oleoperator{\otimes}
 (y \oleoperator{\owedge} z).
$

The repeated ole-addition had been defined in Ref.\ \cite{Borges-2004}, 
and reads
\begin{eqnarray}
 \label{eq:ole-repeated-sum}
  \begin{array}{lll}
  \smallskip
  n \oleoperator{\odot} y &=& 
                               \stackrel[ \scalebox{0.7}{[\textit{q}]} ]
                                        {}{} \!
                               \displaystyle\sum^n y,
  \\
                           &=& 
                               \stackrel[ \scalebox{0.7}{\textit{q}} ]
                                        {}{} \!
                               \bigg[\;
                                     \displaystyle\sum^n \; [y]_q\;
                               \bigg],
 \medskip
  \\
                          &=& \displaystyle
                              \frac{{\big(1+(1-q)y\big)}_+^n -1}
                                   {1-q}.
  \end{array}
\end{eqnarray}
This is identical to Eq.\ (8) of Ref.\ \cite{lemans-2009}.
Analytical extension into the real domain yields the non commutative
ole-dot-multiplication:
\begin{eqnarray}
 \label{eq:ole-dot-product}
  x \oleoperator{\odot} y  &=& \frac{ {\big(1+(1-q)y\big)}_+^x -1 } 
                                    {1-q}.
\end{eqnarray}
The ole-dot-multiplication with the unity is expressed by
$x \oleoperator{\odot} 1 = {}_q\big[\, x \; [1]_q \,\big]$,
with 
$[1]_q = \ln \exp_q(1) = {(1-q)}^{-1} \, \ln(2-q) \ne 1$ for $q \ne 1$ 
and $q < 2$.
This relation connects the ole-dot-multiplication and the le deformed numbers
with the $Q$-analog of $n$ (\ref{eq:quantum-number}):
$n \oleoperator{\odot} 1 = (Q^n - 1)/(Q-1)$, with $Q=2-q$.
The ole-dot power naturally follows from the repeated ole-dot-multiplication,
not shown here.

\subsection{\label{subsec:iel-algebra}iel-Arithmetics}

According to the generating rule for i-algebras 
(\ref{eq:generating-rule-i}), 
the iel-addition is
\begin{eqnarray}
 \label{eq:iel-sum}
 \begin{array}{lll}
 \smallskip
 x \ieloperator{\oplus} y &=&  \left\{ \, 
                                       {}_q\{x\} + {}_q\{y\} \: 
                               \right\}_q,
 \\
                        &=& \displaystyle
                            \text{sign}(x+y) \,
                            \exp \left(\ln_q \left| 
                                             \text{sign}(x)\,\text{e}_q^{\ln|x|}
                                             +
                                             \text{sign}(y)\,\text{e}_q^{\ln|y|}
                                             \right|
                                 \right).
 \end{array}
\end{eqnarray}
The cutoff of the $q$-exponential (\ref{eq:q-exponential})
imposes restrictions on the domain of (\ref{eq:iel-sum}).
Its neutral iel-additive element 
$ N_{\{+\}} = \{0\}_q $, is
\begin{eqnarray}
 \label{eq:iel-neutral-addition}
  \begin{array}{ll}
   N_{\{+\}} \to 0, & q \ge 1,
   \\ \noalign{\smallskip}
   N_{\{+\}} \le \displaystyle \text{e}^{\frac{-1}{1-q}}, & q < 1.
  \end{array}
\end{eqnarray}
For $q < 1$, there are infinite neutral iel-additive elements, 
including the zero.
The iel-difference reads
\begin{eqnarray}
 \label{eq:iel-difference}
 \begin{array}{lll}
 \smallskip
 x \ieloperator{\ominus} y &=&  \left\{ \, 
                                        {}_q\{x\} - {}_q\{y\} \: 
                                \right\}_q,
 \\
                        &=& \displaystyle
                            \text{sign}(x-y) \,
                            \exp \left(\ln_q \left| 
                                             \text{sign}(x)\,\text{e}_q^{\ln|x|}
                                             -
                                             \text{sign}(y)\,\text{e}_q^{\ln|y|}
                                             \right|
                                 \right).
 \end{array}
\end{eqnarray}
The opposite iel-additive element is
\begin{eqnarray}
 \label{eq:iel-opposite}
 \ieloperator{\ominus} y = \left\{
                        \begin{array}{ll}
                            \left.
                             \begin{array}{ll}
                              - y, 
                                  & \text{ if } |y| > \exp\left(
                                                                \frac{-1}{1-q}
                                                          \right)\!,
                              \\ \noalign{\smallskip}
                              - \text{sign}(y)\,\exp\left(
                                                          \frac{-1}{1-q}
                                                    \right)\!,
                                  & \text{ if } |y| \le \exp\left(
                                                                  \frac{-1}{1-q}
                                                            \right)\!,
                             \end{array}
                            \right\}
                            & q < 1,
                        \\ \noalign{\smallskip}
                        \; -y, & q = 1,
                        \\ \noalign{\smallskip}
                            \left.
                             \begin{array}{ll}
                              - y, 
                              & \text{ if } |y| < \exp\left(
                                                            \frac{1}{q-1}
                                                      \right)\!,
                              \\ \noalign{\smallskip}
                              - \text{sign}(y)\, \exp\left(
                                                           \frac{1}{q-1}
                                                     \right)\!,
                              & \text{ if }|y|\ge \exp\left(
                                                            \frac{1}{q-1}
                                                      \right)\!,
                             \end{array}
                            \right\}
                            & q > 1.
                        \end{array}
                           \right.
\end{eqnarray}
The iel-multiplication and the iel-division are
\begin{eqnarray}
 \label{eq:iel-product}
 \begin{array}{lll}
 \smallskip
 x \ieloperator{\otimes} y &=&  \left\{ \, 
                                        {}_q\{x\} \; {}_q\{y\} \: 
                                 \right\}_q,
 \\
                        &=& \displaystyle
                            \text{sign}(xy) \,
                            \exp \left(
                                       \ln_q \left( 
                                                   \text{e}_q^{\ln|x|} 
                                                   \; 
                                                   \text{e}_q^{\ln|y|}
                                             \right) 
                                 \right),
 \end{array}
\end{eqnarray}
\begin{eqnarray}
 \label{eq:iel-ratio}
 \begin{array}{lll}
 \smallskip
 x \ieloperator{\oslash} y &=&  \left\{ \, 
                                        \displaystyle
                                        \frac{{}_q\{x\}}{{}_q\{y\}} \: 
                                \right\}_q,
 \\
                        &=& \displaystyle
                            \text{sign}(x/y) \,
                            \exp \left(
                                       \ln_q \left(
                                              \frac{\text{e}_q^{\ln|x|}}
                                                   {\text{e}_q^{\ln|y|}}
                                             \right) 
                                 \right).
 \end{array}
\end{eqnarray}
The neutral element of the iel-multiplication is 
$I_{\scalebox{.57}{\{$\times$\}}} = \{1\}_q = 1$.
The iel-absorbing element coincides with the neutral iel-additive element,
$A_{\scalebox{.57}{\{$\times$\}}} = N_{\scalebox{.57}{\{+\}}}$
(\ref{eq:iel-neutral-addition}).

The repeated iel-multiplication (iel-power) is given by
\begin{eqnarray}
 \label{eq:iel-repeated-product}
  x {\ieloperator{\owedge} n} =
                                \ieloperator{\displaystyle\prod^n} x
                              = \left\{ \;
                                       \left( 
                                             {}_q\{x\} \;
                                       \right)^n \:
                                \right\}_q,
\end{eqnarray}
that is rewritten as
(after analytical extension from $n \in \mathbb{N}$ to $y \in \mathbb{R}$)
\begin{eqnarray}
 \label{eq:iel-power}
  x {\ieloperator{\owedge} y} = 
                               \exp\big(
                                        \ln_q ( \, \text{e}_q^{\ln |x|} \,)^y \,
                                   \big),
                               \quad (x>0),
\end{eqnarray}
with the particular cases
$x {\ieloperator{\owedge} 0} = 1$ ($x \ne 0$), \ 
$x {\ieloperator{\owedge} 1} = x$ ($x \ne 0$), \ 
$1 {\ieloperator{\owedge} y} = 1$ ($y \ne 0$), \ 
$\lim_{x\to 0^+} (x {\ieloperator{\owedge} y}) = \exp\big(-1/(1-q)\big)$ 
                                                         ($y > 0$ and $q<1$), \ 
$\lim_{x\to 0^+} (x {\ieloperator{\owedge} y}) \to \infty$ 
                                                         ($y < 0$ and $q<1$), \ 
$\lim_{x\to 0^+} (x {\ieloperator{\owedge} y}) = 0$        
                                                         ($y > 0$ and $q>1$), \ 
$\lim_{x\to 0^+} (x {\ieloperator{\owedge} y}) = \exp\big(-1/(1-q)\big)$ 
                                                         ($y < 0$ and $q>1$), \ 
$x\,{\owedge_{\{1\}}\, y} = x^y$.
The iel-power is right-distributive with respect to the iel-multiplication:
$
 (x \ieloperator{\otimes} y) \ieloperator{\owedge} z 
 =
 (x \ieloperator{\owedge} z)
 \ieloperator{\otimes}
 (y \ieloperator{\owedge} z).
$

The repeated iel-addition defines the iel-dot-multiplication:
\begin{eqnarray}
 \label{eq:iel-repeated-sum}
  \begin{array}{lll}
  \smallskip
  n \ieloperator{\odot} y &=& \ieloperator{\displaystyle\sum^n} y,
  \\
                          &=& \bigg\{\:
                                     \displaystyle\sum^n \; {}_q\{y\}\;
                              \bigg\}_q,
  \smallskip
  \\
                          &=& \text{sign}(y)\, \exp(\ln_q n) \: |y|^{n^{1-q}}.
  \end{array}
\end{eqnarray}
Analytical extension from $n \in \mathbb{N}$ to $x \in \mathbb{R}_+$
can be represented by 
\begin{equation}
 \label{eq:iel-dot-product}
  x \ieloperator{\odot} y  = \text{sign}(y)\, 
                             \exp (\ln_q x) \, |y|^{x^{1-q}},
                             \quad (x>0).
\end{equation}
The iel-number is connected to the iel-dot-multiplication by
$x \ieloperator{\odot} 1 =  \big\{\, x \; {}_q\{1\}\, \big\}_q = \{x\}_q$,
since ${}_q\{1\}=1$.

\subsection{\label{subsec:oel-algebra}oel-Arithmetics}

The oel-arithmetic operators derives from (\ref{eq:generating-rule-o}):
\begin{eqnarray}
 \label{eq:oel-sum}
 x \oeloperator{\oplus} y &=&  {}_q\!\left\{ \, 
                                             \{x\}_{q} + \{y\}_{q} \, 
                                     \right\},
 \nonumber
 \\
                        &=& \displaystyle
                            \text{sign}(x+y) \,
                            \exp_q \left(\ln \left| 
                                             \text{sign}(x)\,\text{e}^{\ln_q|x|}
                                             +
                                             \text{sign}(y)\,\text{e}^{\ln_q|y|}
                                             \right|
                                 \right),
\\
 \label{eq:oel-difference}
 x \oeloperator{\ominus} y &=&  {}_q\!\left\{ \, 
                                              \{x\}_{q} - \{y\}_{q} \, 
                                      \right\},
 \nonumber
 \\
                        &=& \displaystyle
                            \text{sign}(x-y) \,
                            \exp_q \left(\ln \left| 
                                             \text{sign}(x)\,\text{e}^{\ln_q|x|}
                                             -
                                             \text{sign}(y)\,\text{e}^{\ln_q|y|}
                                             \right|
                                 \right),
\\
 \label{eq:oel-product-generating-rule}
 x \oeloperator{\otimes} y &=&  {}_q\!\left\{ \, 
                                              \{x\}_{q} \: \{y\}_{q}  
                                      \right\},
 \nonumber
 \\
                        &=& \displaystyle
                            \text{sign}(xy) \,
                            \exp_q \left(\ln \left| 
                                                   \text{e}^{\ln_q|x|}  \:
                                                   \text{e}^{\ln_q|y|}
                                             \right|
                                 \right),
\\
 \label{eq:oel-ratio-generating-rule}
 x \oeloperator{\oslash} y &=&  {}_q\!\left\{ \, \displaystyle
                                             \frac{ \{x\}_{q} }{  \{y\}_{q} }\,
                                      \right\},
 \nonumber
 \\
                        &=& \displaystyle
                            \text{sign}(x/y) \,
                            \exp_q \left(\ln \left| 
                                                   \frac{\text{e}^{\ln_q|x|}}
                                                        {\text{e}^{\ln_q|y|}}
                                             \right|
                                 \right).
\end{eqnarray}

Equations (\ref{eq:oel-product-generating-rule}) 
      and (\ref{eq:oel-ratio-generating-rule}) 
can be rearranged as
\begin{eqnarray}
 \label{eq:oel-product}
 x \oeloperator{\otimes} y &=& \displaystyle 
                               \text{sign}(xy) \,
                               \left( 
                                     |x|^{1-q} + |y|^{1-q} - 1 
                               \right)_+^{\frac{1}{1-q}}
\end{eqnarray}
and
\begin{eqnarray}
 \label{eq:oel-ratio}
 x \oeloperator{\oslash} y &=& \displaystyle
                               \text{sign}(x/y) \,
                               \left( 
                                     |x|^{1-q} - |y|^{1-q} + 1 
                               \right)_+^{\frac{1}{1-q}}.
\end{eqnarray}

\begin{figure}[!htb]
\centering
 \begin{minipage}[h]{0.33\linewidth}
 \end{minipage}
 \begin{minipage}[h]{0.33\linewidth}
  \includegraphics[width=\linewidth,keepaspectratio,clip=]{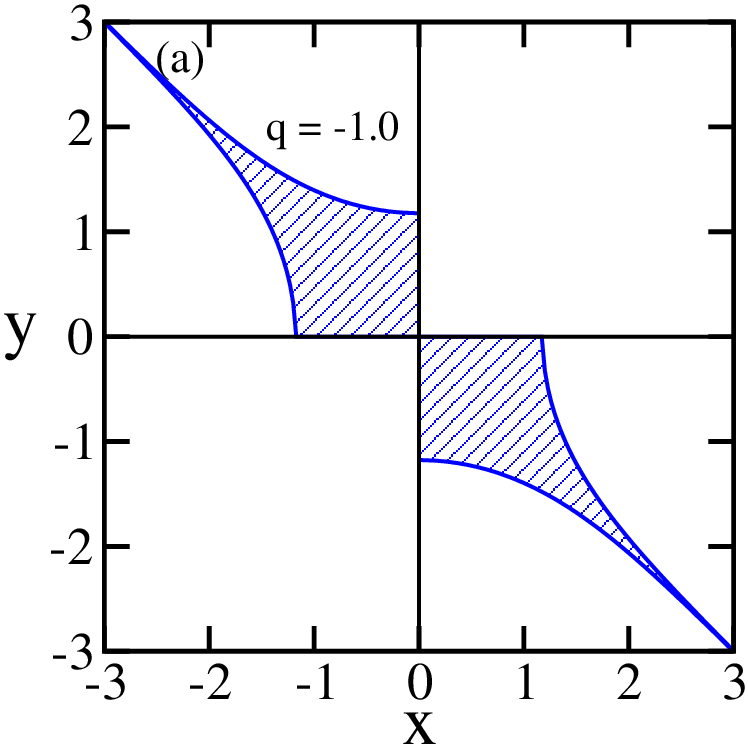}
 \end{minipage}
\qquad \qquad
 \begin{minipage}[h]{0.33\linewidth}
  \includegraphics[width=\linewidth,keepaspectratio,clip=]{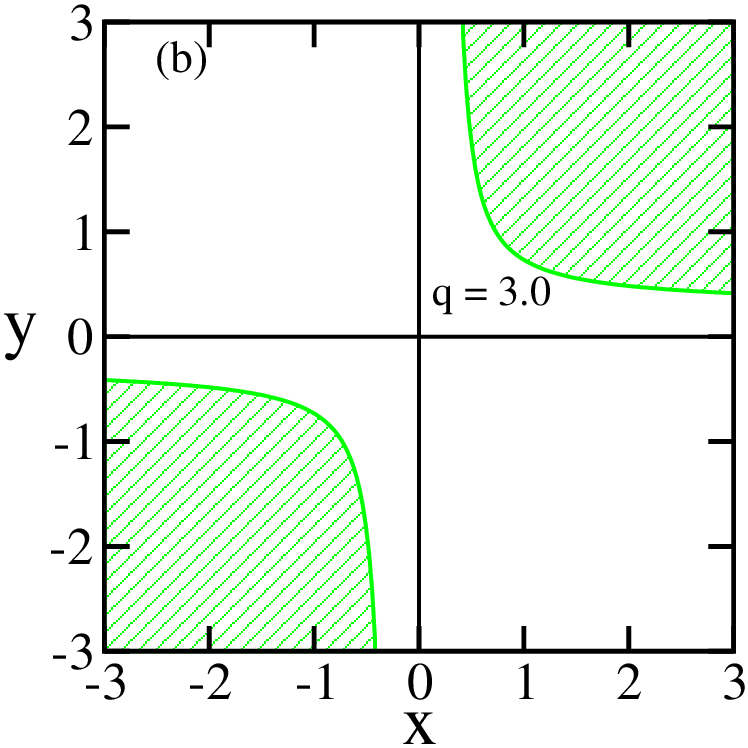}
 \end{minipage}
\\
\medskip
 \begin{minipage}[h]{0.33\linewidth}
 \end{minipage}
 \begin{minipage}[h]{0.33\linewidth}
  \includegraphics[width=\linewidth,keepaspectratio,clip=]{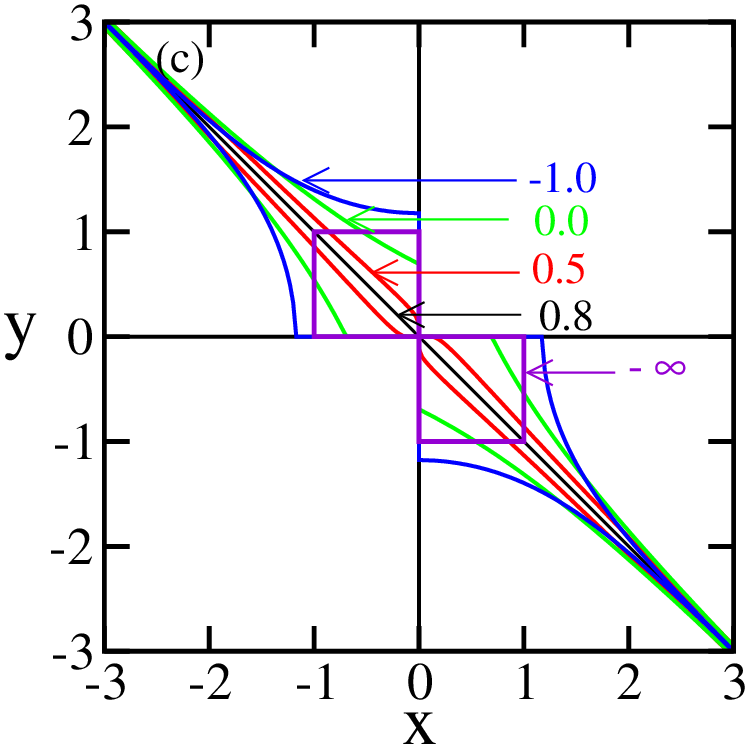}
 \end{minipage}
\qquad \qquad
 \begin{minipage}[h]{0.33\linewidth}
  \includegraphics[width=\linewidth,keepaspectratio,clip=]{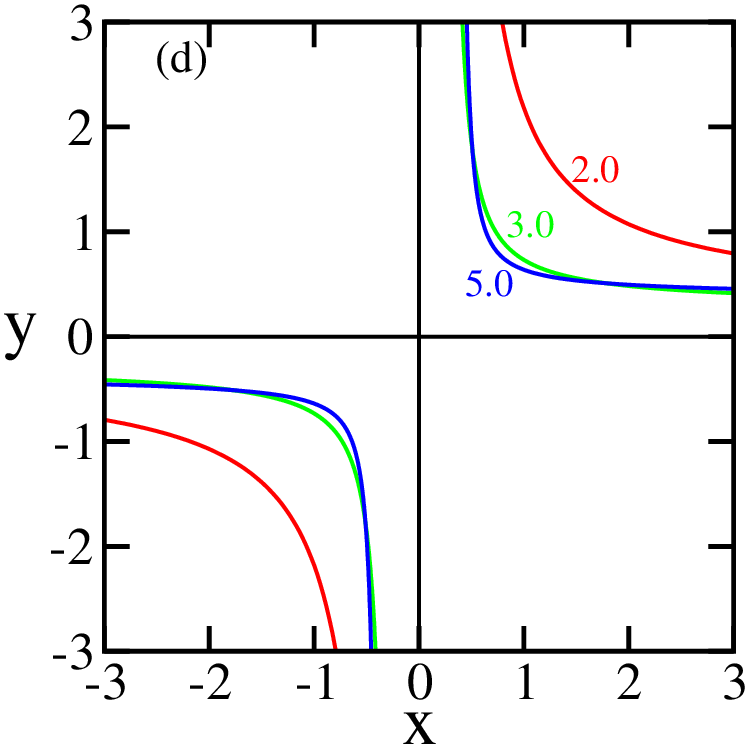}
 \end{minipage}
 \caption{
          Cutoff regions for the oel-addition 
          (\protect\ref{eq:oel-sum}).
          Left column: $q<1$, right column: $q>1$.
          Top line: 
          the shaded regions correspond to the cutoff regions 
          of the oel-addition.
          (a) $q=-1$.
          (b) $q=3$.
          Bottom line: the curves represent the cutoff borders.
          Regions are not shaded to avoid excessively heavy representation.
          Their pattern is similar to (a) or (b):
          for $q<1$, the cutoff regions lie inside the corresponding 
          closed curves,
          and for $q>1$, the cutoff regions lie outside the corresponding 
          curves.
          (c) Different values of $q<1$ (indicated). 
              The cutoff region shrinks and eventually collapses at $y=-x$
              as $q \to 1^{-}$.
          (d) Different values of $q>1$ (indicated). 
              As $q \to 1^{+}$, the non connected regions depart from the 
              origin, and there are no cutoff regions.
          \label{fig:oel-sum-cutoff}
         }
\end{figure}

\begin{figure}[!htb]
\centering
 \begin{minipage}[h]{0.33\linewidth}
 \end{minipage}
 \begin{minipage}[h]{0.33\linewidth}
  \includegraphics[width=\linewidth,keepaspectratio,clip=]{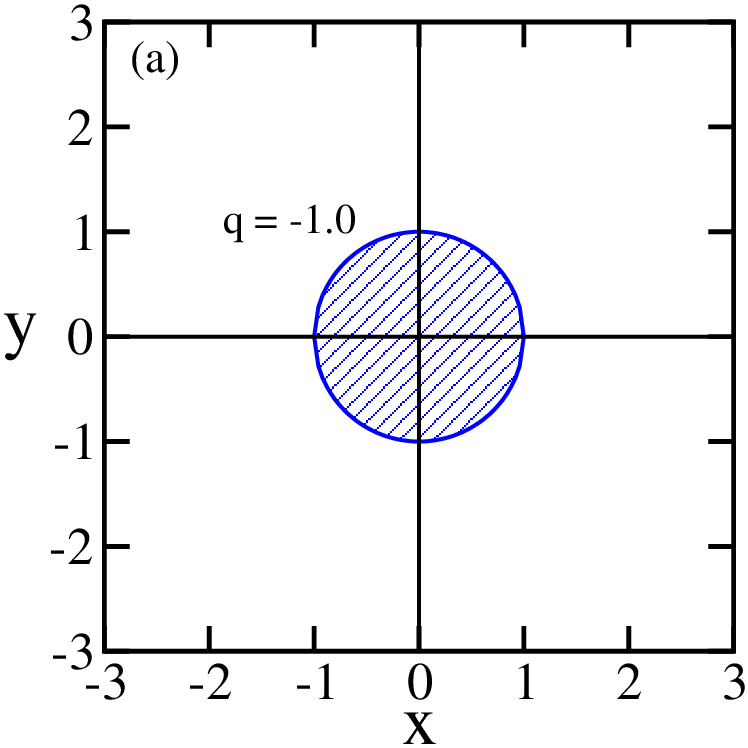}
 \end{minipage}
\qquad \qquad
 \begin{minipage}[h]{0.33\linewidth}
  \includegraphics[width=\linewidth,keepaspectratio,clip=]{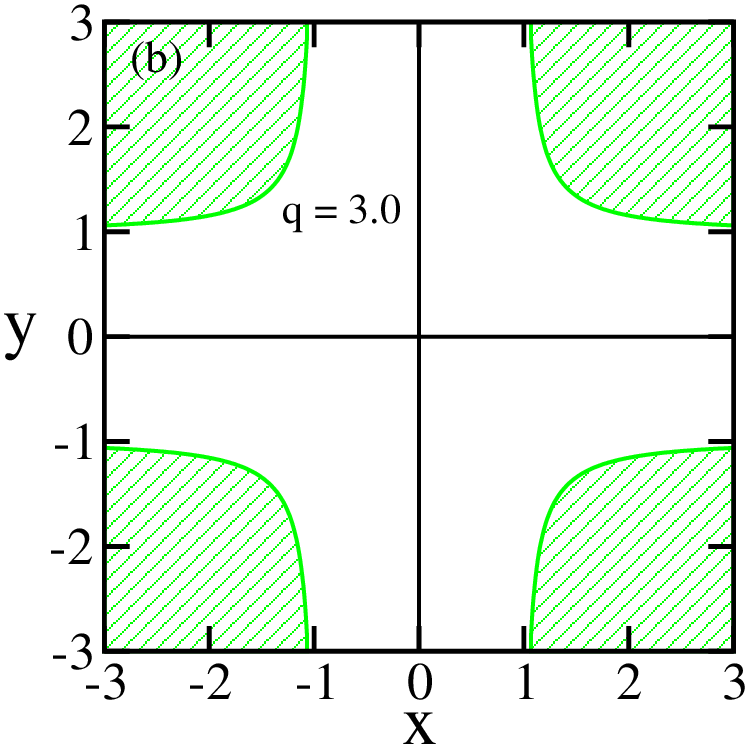}
 \end{minipage}
\\
\medskip
 \begin{minipage}[h]{0.33\linewidth}
 \end{minipage}
 \begin{minipage}[h]{0.33\linewidth}
  \includegraphics[width=\linewidth,keepaspectratio,clip=]{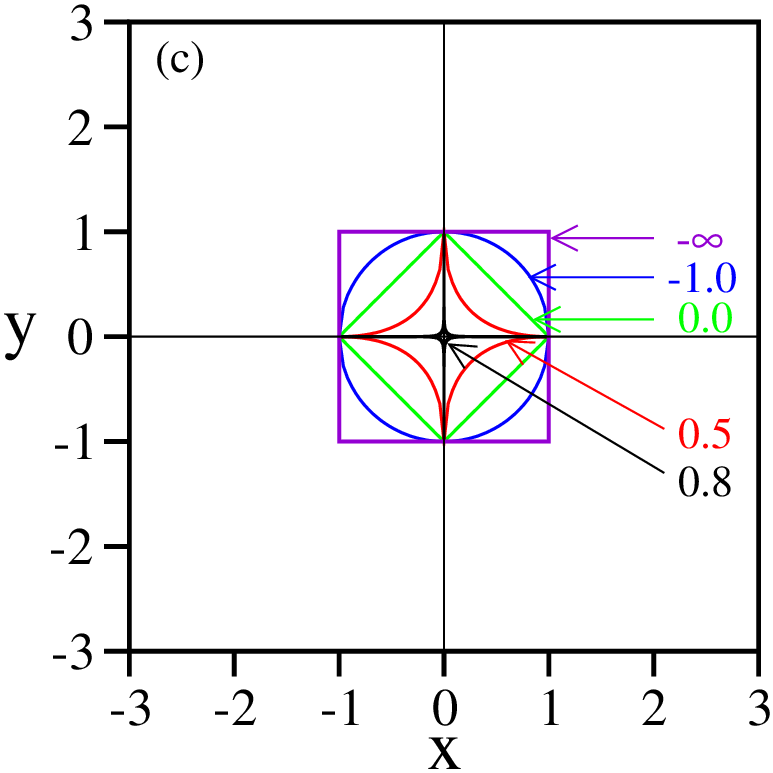}
 \end{minipage}
\qquad \qquad
 \begin{minipage}[h]{0.33\linewidth}
  \includegraphics[width=\linewidth,keepaspectratio,clip=]{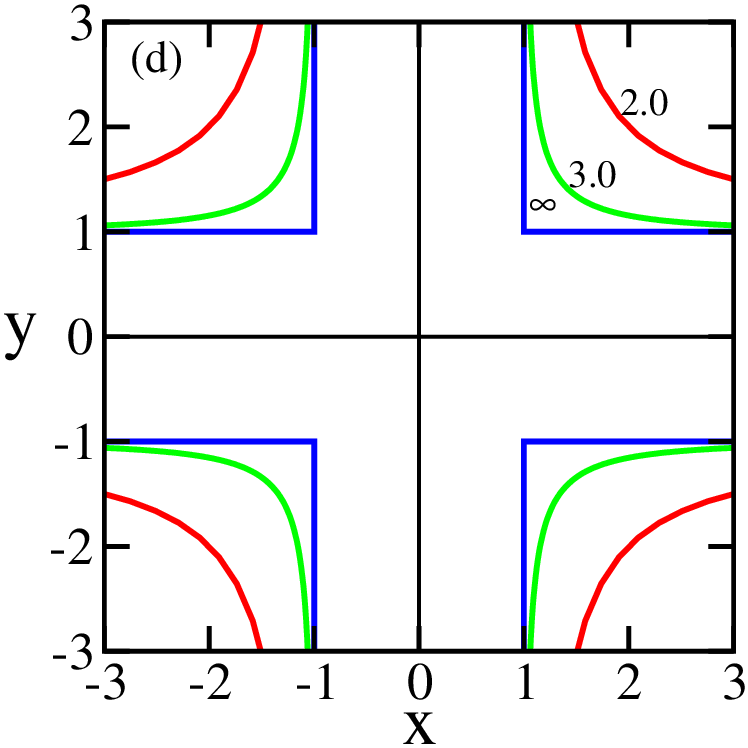}
 \end{minipage}
 \caption{
          Cutoff regions for the oel-multiplication 
          (\protect\ref{eq:oel-product-generating-rule}).
          Left column: $q<1$, right column: $q>1$.
          Top line: 
          the shaded regions correspond to the cutoff regions 
          of the oel-multiplication.
          (a) $q=-1$.
          (b) $q=3$.
          Bottom line: the curves represent the cutoff borders,
          $|y| = (1 - |x|^{1-q})^{1/(1-q)}$.
          Regions are not shaded to avoid excessively heavy representation.
          Their pattern is similar to the adopted in (a) or (b):
          for $q<1$, the cutoff regions lie inside the corresponding 
          closed curves,
          and for $q>1$, the cutoff regions lie outside the corresponding 
          curves.
          (c) Different values of $q<1$ (indicated). 
              The cutoff region shrinks and eventually collapses at $(0,0)$ 
              as $q \to 1^{-}$,
              when the curves coincide with the axes.
          (d) Different values of $q>1$ (indicated). 
              As $q \to 1^{+}$, the non connected regions depart from the 
              origin, and there are no cutoff regions.
          \label{fig:oel-product-cutoff}
         }
\end{figure}

The oel-product and the oel-ratio were defined in Ref.\ \cite{Borges-2004},
referred to as $q$-product and $q$-ratio, respectively
(see also Subsec.\ 3.3.2 of Ref.\ \cite{Tsallis-springer-2009}).
The cutoff that appears in (\ref{eq:oel-product}) defines 
regions in which the oel-arithmetical operators are ill-defined.
Figures \ref{fig:oel-sum-cutoff} and \ref{fig:oel-product-cutoff}
show the regions for which the cutoff applies for the oel-addition
and oel-multiplication, respectively.
The first column of each (figures a and c) shows instances for $q<1$, 
and the second column    (figures b and d), for $q>1$.
The first line (figures a and b) exhibits the cutoff regions
with a shaded pattern for one typical value of the parameter $q$. 
The second line (figures c and d) display superimposed curves 
of the borders of the cutoff regions for various values of $q$, 
without shading them, otherwise they would be confusing;
they follow the same pattern of the corresponding figures a and b, 
respectively.
The cutoff regions are closed for $q<1$ 
(illustrated with $q=-1$ by
 figures \ref{fig:oel-sum-cutoff}a and \ref{fig:oel-product-cutoff}a),
and they are open and not connected, 
lying on the outer side delimited by the bounding curves, for $q>1$
(illustrated with $q=3$ by
 figures \ref{fig:oel-sum-cutoff}b and \ref{fig:oel-product-cutoff}b).
The second line of the figures help us to understand the effect of the
deforming parameter $q$ on the cutoff regions.
As $q$ approaches unity from below
 (figures \ref{fig:oel-sum-cutoff}c and \ref{fig:oel-product-cutoff}c),
the cutoff regions become smaller and eventually vanish.
For the oel-addition, Figure \ref{fig:oel-sum-cutoff}c, 
the borders of the cutoff region approach the second bisector ($y=-x$),
and, 
for the oel-multiplication, Figure \ref{fig:oel-product-cutoff}c, 
they approach the origin $(0,0)$.
As $q$ approaches unity from above 
 (figures \ref{fig:oel-sum-cutoff}d and \ref{fig:oel-product-cutoff}d),
the cutoff regions move away from the origin. 
At $q=1$, no pair of numbers $(x,y)$ fall within the cutoff regions,
and the ordinary arithmetic operators are defined everywhere.
%
%

%

The distributivity of the oel-multiplication with respect to the oel-addition
is valid whenever the cutoff conditions of the l.h.s.\ and the r.h.s.\ of 
$x \oeloperator{\otimes} (y \oeloperator{\oplus} z)
 =
 (x \oeloperator{\otimes} y) \oeloperator{\oplus} (x \oeloperator{\otimes} z) $
are not met.
As $q$ approaches unity, even from below or from above, the distributivity
of the oel-multiplication with respect to the oel-addition is valid 
for all real values $(x,y,z)$.

The neutral oel-additive element 
is $ {}_{\scalebox{.57}{\{+\}}}\!N = {}_q\{0\} = 0 $ for $q \ge 1$,
and $\nexists \: {}_{\scalebox{.57}{\{+\}}}\!N \:| \:
{}_{\scalebox{.57}{\{+\}}}\!N \oeloperator{\oplus} \: x = x$ for $q < 1$.
As a consequence, there is no opposite oel-additive element for $q < 1$.
For $q \ge 1$, $\oeloperator{\ominus} y = -y$.
The absorbing element 
$ {}_{\scalebox{.57}{\{$\times$\}}}\!A = {}_q\{0\} = 0 $ for $q \ge 1$,
and
$\nexists \: {}_{\scalebox{.57}{\{$\times$\}}}\!A \:| \:
{}_{\scalebox{.57}{\{$\times$\}}}\!A \oeloperator{\otimes} \: x = 0$ for $q<1$
and $|x| > 1$. 
If $q<1$, and $|x| < 1$ 
the cutoff of (\ref{eq:oel-product})
[see (\ref{eq:q-exponential})]
implies that zero is an absorbing element,
and, in this case,
differently from the other three generalized algebras,
${}_{\scalebox{.57}{\{$+$\}}}\!N \ne {}_{\scalebox{.57}{\{$\times$\}}}\!A$.
The neutral multiplicative element of the oel-multiplication is
${}_{\scalebox{0.57}{\{$\times$\}}}I = {}_q\{1\} = 1$, for all values of $q$.
The inverse oel-multiplicative element is
\begin{eqnarray}
 \label{eq:oel-inverse-multiplicative}
 1 \! \oeloperator{\oslash} y &=& \left\{
                                     \begin{array}{ll}
                                      \text{sign}(y)\, \left( 
                                                             2 - |y|^{1-q}
                                                       \right)^{\frac{1}{1-q}},
                                         & \text{if }|y| < 2^{\frac{1}{1-q}},
                                      \\ \noalign{\smallskip}
                                      0, & \text{otherwise}.
                                     \end{array}
                               \right. 
\end{eqnarray}
This implies the unorthodox property
$\lim_{y\to 0^+} (1 \!\oeloperator{\oslash} y) \to 2^{1/(1-q)}$, for $q < 1$.

The oel-power, 
previously defined in Ref.\ \cite{Borges-2004} (with different symbols), 
is written as
\begin{eqnarray}
 \label{eq:oel-repeated-product}
  x {\oeloperator{\owedge} n} =
                                \stackrel[ \scalebox{0.7}{\{\textit{q}\}} ]
                                         {}{} \!\!
                                \displaystyle\prod^n x
                              = {}_q\left\{ \:
                                           \left(
                                                 \{x\}_q
                                           \right)^n
                                    \right\},
                                  \quad (x>0).
\end{eqnarray}
This operator also appears as Eq.\ (8) of Ref.\ \cite{lemans-2009}. 
We make an analytical extension from $n \in \mathbb{N}$ to $y \in \mathbb{R}$,
and the oel-power can also be written as
\begin{eqnarray}
 \label{eq:oel-power}
  x {\oeloperator{\owedge} y} = 
                                  \exp_q\left(\, 
                                              y \ln_q x \,
                                        \right),
                                  \quad (x>0).
\end{eqnarray}
Particular cases are 
$x {\oeloperator{\owedge} 0} = 1$ ($x \ne 0$), \ 
$x {\oeloperator{\owedge} 1} = x$ ($x \ne 0$), \ 
$1 {\oeloperator{\owedge} y} = 1$ ($y \ne 0$), \ 
$\lim_{x\to 0^+} (x{\oeloperator{\owedge} y}) = 0$        ($y \ge 1$, $q<1$), \ 
$\lim_{x\to 0^+} (x{\oeloperator{\owedge} y}) = \exp_q\big(-y/(1-q)\big)$ 
                                                          ($y < 1$,   $q<1$), \ 
$\lim_{x\to 0^+} (x{\oeloperator{\owedge} y}) = 0$        ($y > 0$,   $q>1$), \ 
$\lim_{x\to 0^+} (x{\oeloperator{\owedge} y}) \to \infty$ ($y < 0$,   $q>1$), \ 
and, as always,
$x\,{{}_{\{1\}}{\owedge} \; y} = x^y$.
The oel-power is right-distributive with respect to the oel-multiplication:
$
 (x \oeloperator{\otimes} y) \oeloperator{\owedge} z 
 =
 (x \oeloperator{\owedge} z)
 \oeloperator{\otimes}
 (y \oeloperator{\owedge} z).
$

The repeated oel-addition is
\begin{equation}
 \label{eq:oel-repeated-sum}
 \stackrel[ \scalebox{0.7}{\{\textit{q}\}} ] {}{} \!\!
 \displaystyle\sum^n y 
                      = \stackrel[ \scalebox{0.7}{\textit{q}} ]
                                 {}{} \!
                        \left\{\;\displaystyle\sum^n \; \{y\}_q\;\right\}.
\end{equation}
Its analytical extension from $n \in \mathbb{N}$ to $x \in \mathbb{R}_+$ 
defines the non commutative oel-dot-multiplication:
\begin{eqnarray}
 \label{eq:oel-dot-product}
  \begin{array}{lll}
  x \oeloperator{\odot} y = \text{sign}(y)\, \left( 
                                                     (1-q) \ln x + |y|^{1-q} 
                                              \right)_+^{\frac{1}{1-q}},
                            \quad (x>0).
  \end{array}
\end{eqnarray}
The oel-number is connected to the oel-dot-multiplication by
$x \oeloperator{\odot} 1 = {}_q\big\{\, x \; \{1\}_q \, \} = {}_q\{x\}$,
since $\{1\}_q=1$.

\section{\label{sec:q-calculus}Deformed $q$-calculus}

Following the lines of Ref.\ \cite{Borges-2004}
(see also sections II.C and II.D of \cite{daCosta-Gomez-Borges-2020}),
we connect the deformed algebra with deformed calculus,
and define the deformed differentials of ordinary numbers:
\begin{subequations}
 \label{eq:q-differentials-def}
 \begin{equation}
 \label{eq:ile-differential-def}
  \ileoperator{\text{d}}\! x = \lim_{x' \to x} \big( 
                                                    x' \ileoperator{\ominus} x
                                                \big),
 \end{equation}
 \begin{equation}
 \label{eq:ole-differential-def}
  \oleoperator{\text{d}}\! x = \lim_{x' \to x} \big( 
                                                    x' \oleoperator{\ominus} x
                                                \big),
 \end{equation}
 \begin{equation}
 \label{eq:iel-differential-def}
  \ieloperator{\text{d}}\! x = \lim_{x' \to x} \big( 
                                                    x' \ieloperator{\ominus} x
                                                \big),
 \end{equation}
 \begin{equation}
 \label{eq:oel-differential-def}
  \oeloperator{\text{d}}\! x = \lim_{x' \to x} \big( 
                                                    x' \oeloperator{\ominus} x
                                                \big).
 \end{equation}
\end{subequations}
The definitions of the corresponding deformed differences,
Eq.\ (\ref{eq:ile-difference}), 
     (\ref{eq:ole-difference}), 
     (\ref{eq:iel-difference}), 
     (\ref{eq:oel-difference}),
lead to
\begin{subequations}
\label{eq:q-differentials-equivalence}
 \begin{equation}
  \label{eq:ile-differential-equivalence}
  \ileoperator{\text{d}}\! x = \text{d}\,\big({}_q[x]\,\big),
 \end{equation}
 \begin{equation}
  \label{eq:ole-differential-equivalence}
  \oleoperator{\text{d}}\! x = \text{d}\,\big(\,[x]_q\big),
 \end{equation}
 \begin{equation}
  \label{eq:iel-differential-equivalence}
  \ieloperator{\text{d}}\! x = \text{d}\,\big({}_q\{x\}\,\big),
 \end{equation}
 \begin{equation}
  \label{eq:oel-differential-equivalence}
  \oeloperator{\text{d}}\! x = \text{d}\,\big(\,\{x\}_q\big),
 \end{equation}
\end{subequations}%
i.e., 
the deformed differential of an ordinary variable 
[l.h.s.\ of (\ref{eq:q-differentials-equivalence})]
is equal to 
the ordinary differential of the corresponding complementary deformed variable
[r.h.s.\ of (\ref{eq:q-differentials-equivalence})]:
the i-differential of a variable is equal to the ordinary differential 
of an o-variable,
          (\ref{eq:ile-differential-equivalence}) 
and       (\ref{eq:iel-differential-equivalence}),
and 
the o-differential of a variable is equal to the ordinary differential 
of an i-variable,
          (\ref{eq:ole-differential-equivalence}) 
and       (\ref{eq:oel-differential-equivalence}).
All the deformed differentials given by (\ref{eq:q-differentials-equivalence})
can be arranged as the product of the ordinary differential $\text{d}x$
by a deforming function $h_{\delta}(x)$,
with $\delta \in \{\text{ile}, \text{ole}, \text{iel}, \text{oel} \}$
representing the deformation
(%
$\ileoperator{\text{d}}\!x = h_{\text{ile}}(x)\, \text{d}\,x$,
$\oleoperator{\text{d}}\!x = h_{\text{ole}}(x)\, \text{d}\,x$,
$\ieloperator{\text{d}}\!x = h_{\text{iel}}(x)\, \text{d}\,x$,
$\oeloperator{\text{d}}\!x = h_{\text{oel}}(x)\, \text{d}\,x$%
).
Their explicit forms are
\begin{subequations}
\label{eq:deforming-functions}
 \begin{equation}
  \label{eq:lambda_ile}
   h_{\text{ile}}(x) = \text{e}^{(1-q)x},
 \end{equation}
 \begin{equation}
  \label{eq:lambda_ole}
  h_{\text{ole}} = \frac{1}{1+(1-q)x},
   \quad \left(x \ne \frac{-1}{1-q}\right),
 \end{equation}
 \begin{equation}
  \label{eq:lambda_iel}
  h_{\text{iel}}(x) = \frac{1}{x} 
                                    \big(1+(1-q)\ln x\big)^{\frac{q}{1-q}},
   \quad (x > 0),
 \end{equation}
 \begin{equation}
  \label{eq:lambda_oel}
  h_{\text{oel}}(x) = \frac{1}{x^q} 
                                    \exp\left(\frac{x^{1-q}-1}{1-q}\right),
   \quad (x > 0).
 \end{equation}
\end{subequations}%

A pair of generalized derivatives of a function $f(x)$,
holding a duality nature between them,
stem from each of the deformed differentials,
according to which variable the deformed differential applies on:
whether on the independent variable $x$,
--- and thus a linear deformed derivative ---, 
generically represented by $\text{D}_\delta f(x)$,
or on the dependent variable $f$,
--- and thus a nonlinear deformed derivative ---,
generically represented by $\widetilde{\text{D}}_\delta f(x)$,
resulting eight different cases:
\\
{\bf 1. ile-derivatives}
\begin{subequations}
\label{eq:ile-derivatives}
\\
linear~ile-derivative:
 \begin{eqnarray}
 \label{eq:linear-ile-derivative}
 \displaystyle
  \text{D}_{\text{ile}} f(x) \equiv
  \frac{\text{d}f(x)}{\ileoperator{\text{d}}x} =
   \frac{1}{h_{\text{ile}}(x)}
   \frac{\text{d}f(x)}{\text{d}x},
 \end{eqnarray}
\\
  \text{{nonlinear~ile-derivative:}} 
 \begin{eqnarray}
 \label{eq:nonlinear-ile-derivative}
 \displaystyle
  \widetilde{\text{D}}_{\text{ile}} f(x) \equiv
  \frac{ \ileoperator{\text{d}} f(x)} {\text{d}x} =
   h_{\text{ile}}\big(f(x)\big) \,
   \frac{\text{d}f(x)}{\text{d}x}.
 \end{eqnarray}
\end{subequations}%
{\bf 2. ole-derivatives}
\begin{subequations}
\label{eq:ole-derivatives}
\\
  \text{{linear~ole-derivative:}} 
 \begin{eqnarray}
 \label{eq:linear-ole-derivative}
 \displaystyle
  \text{D}_{\text{ole}} f(x) \equiv
  \frac{\text{d}f(x)}{\oleoperator{\text{d}}x} =
   \frac{1}{h_{\text{ole}}(x)}
   \frac{\text{d}f(x)}{\text{d}x},
 \end{eqnarray}
\\
  \text{{nonlinear~ole-derivative:}} 
 \begin{eqnarray}
 \label{eq:nonlinear-ole-derivative}
 \displaystyle
  \widetilde{\text{D}}_{\text{ole}} f(x) \equiv
  \frac{ \oleoperator{\text{d}} f(x)} {\text{d}x} =
   h_{\text{ole}}\big(f(x)\big) \,
   \frac{\text{d}f(x)}{\text{d}x}.
 \end{eqnarray}
\end{subequations}%
{\bf 3. iel-derivatives}
\begin{subequations}
\label{eq:iel-derivatives}
\\
  \text{{linear~iel-derivative:}} 
 \begin{eqnarray}
 \label{eq:linear-iel-derivative}
 \displaystyle
  \text{D}_{\text{iel}} f(x) \equiv
  \frac{ \text{d} f(x)} {\ieloperator{\text{d}}x} =
   \frac{1}{h_{\text{iel}}(x)}
   \frac{\text{d}f(x)}{\text{d}x},
 \end{eqnarray}
\\
  \text{{nonlinear~iel-derivative:}} 
 \begin{eqnarray}
 \label{eq:nonlinear-iel-derivative}
 \displaystyle
  \widetilde{\text{D}}_{\text{iel}} f(x) \equiv
  \frac{ \ieloperator{\text{d}} f(x)} {\text{d}x} =
   h_{\text{iel}}\big(f(x)\big) \,
   \frac{\text{d}f(x)}{\text{d}x}.
 \end{eqnarray}
\end{subequations}%
{\bf 4. oel-derivatives}
\begin{subequations}
\\
  \text{{linear~oel-derivative:}} 
\label{eq:oel-derivatives}
  \text{{linear~oel-derivative:}} 
 \begin{eqnarray}
 \label{eq:linear-oel-derivative}
 \displaystyle
  \text{D}_{\text{oel}} f(x) \equiv
  \frac{\text{d}f(x)}{ \oeloperator{\text{d}}x} =
   \frac{1}{h_{\text{oel}}(x)}
   \frac{\text{d}f(x)}{\text{d}x},
 \end{eqnarray}
\\
  \text{{nonlinear~oel-derivative:}} 
 \begin{eqnarray}
 \label{eq:nonlinear-oel-derivative}
 \displaystyle
  \widetilde{\text{D}}_{\text{oel}} f(x) \equiv
  \frac{ \oeloperator{\text{d}} f(x)} {\text{d}x} =
   h_{\text{oel}}\big(f(x)\big) \,
   \frac{\text{d}f(x)}{\text{d}x}.
 \end{eqnarray}
\end{subequations}%

The duality between the linear and the nonlinear generalized derivatives 
is expressed by
$\text{D}_{\delta} f(x) = \widetilde{\text{D}}_{\delta} f^{-1}(x)$. 
The el-derivatives are defined for $x>0$.
The ole-derivatives had been defined in Ref.\ \cite{Borges-2004},
then referred to as $q$-derivative (the linear deformed derivative)
and its dual $q$-derivative     (the nonlinear deformed derivative).
Particularly, the linear ole-derivative 
(\ref{eq:linear-ole-derivative})
was used to generalize Fisher's information measure 
and the Cramer-Rao inequality \cite{Pennini-Plastino-Ferri-2008}.
The eigenfunction of the linear i/o-deformed derivative
is the ordinary exponential of the o/i-deformed variable,
what directly follows from (\ref{eq:q-differentials-equivalence}).
They are
(written with the symbols $\langle \cdot \rangle$ 
 representing either $[\cdot]$ or $\{\cdot\}$)
\begin{subequations}
 \label{eq:eigenfunctions}
 \begin{eqnarray}
  \label{eq:i-eigenfunctions}
  \frac{\text{d} \, \exp\left({}_q\langle x \rangle \right) }
       {\text{d}_{\langle q \rangle} \: x}   
 =
  \frac{\text{d} \, \exp\left({}_q\langle x \rangle \right) }
       {\text{d} \: \big({}_q{\langle x \rangle}\big)} 
 =
  \exp\left(\langle x \rangle_q \right) 
 \end{eqnarray}
and
 \begin{eqnarray}
  \label{eq:o-eigenfunctions}
  \frac{  \text{d} \, \exp\left( \langle x \rangle_q \right)  }
       {{}_{\langle q \rangle}\text{d} \: x} 
 =
  \frac{  \text{d} \, \exp\left( \langle x \rangle_q \right)  }
       {\text{d} \: \big({\langle x \rangle}_q\big)}
 = 
  \exp\left( \langle x \rangle_q \right).
 \end{eqnarray}
\end{subequations}%
Particularly, the $q$-exponential (\ref{eq:q-exponential})
is the eigenfunction of the linear ole-derivative,
$\text{D}_{\text{ole}} \text{e}_q^x = \text{e}_q^x$
[a particular case of (\ref{eq:o-eigenfunctions}) 
 with $\text{e}_q^x = \text{e}^{[x]_q} $, see (\ref{eq:exp-ilex})].
Alternatively, its ordinary derivative is
${\text{d}\,\text{e}_q^x}\,/\,{\text{d}\,x} = \big(\text{e}_q^x\big)^q$.
The nonlinear deformed derivative of which the $q$-exponential is eigenfunction
was defined in Ref.\ \cite{Nobre-RegoMonteiro-Tsallis-2011}:
\begin{equation}
 \label{eq:nobre-nonlinear-derivative}
 \widetilde{\mathfrak{D}}_{q} f(u) = [f(u)]^{1-q} \, \frac{df(u)}{du},
\end{equation}
where we have used the symbol, $\widetilde{\mathfrak{D}}_{q}$
to distinguish it from the present deformed derivatives.

The integral of the inverse of a variable, $\int_1^x t^{-1} \text{d}t$,
is typically associated to, and frequently taken as the definition of,
the logarithm function.
The general nonlinear cases are
 \begin{eqnarray}
 \label{eq:nonlinear-derivative-logarithms}
  \frac{\text{d}_{\langle q \rangle} \; \langle \ln x \rangle_q}
       {\text{d}\,x}
 =
  \frac{{}_{\langle q \rangle}{\text{d}} \; {}_q\langle \ln x \rangle}
       {\text{d}\,x}
 = \frac{1}{x}.
 \end{eqnarray}
The particular case of this equation for the nonlinear ole-derivative is
[see (\ref{eq:qlog-as-log})]:
$\widetilde{\text{D}}_{\text{ole}} \ln_q x = 1/x$.
Alternatively, the ordinary derivative of the $q$-logarithm is
${\text{d}\ln_q x}\,/\,{\text{d}\,x} = {1}/{x^q}$.
This expression yields an integral representation 
of the $q$-logarithm function,
\begin{equation}
 \label{eq:qlog-integral-representation}
 \int_1^x t^{-q} \text{d}t = \ln_q x.
\end{equation}
The dual linear deformed derivative of (\ref{eq:nobre-nonlinear-derivative}),
defined by Eq.\ (25) of Ref.\ \cite{daCosta-Gomez-Borges-2020},
\begin{equation}
 \label{eq:nobre-linear-derivative}
 \mathfrak{D}_{q} f(x) = \frac{1}{x^{1-q}} \, \frac{\text{d}f(x)}{\text{d}x},
\end{equation}
operates on the $q$-logarithm similarly to the nonlinear ole-derivative:
$\mathfrak{D}_{q} \ln_q x = 1/x$.

Generalized derivatives of a power (for the linear case), 
or generalized powers (for the nonlinear case), of $q$-numbers, are
\begin{eqnarray}
 \label{eq:D-power-linear}
 \begin{array}{llll}
\smallskip
 &\text{D}_{\text{i}} \, \big(\,{}_q\langle x \rangle^n \,\big)
  &=& n \; {}_q\langle x \rangle^{n-1},
\\
 &\text{D}_{\text{o}} \,\big(\,\langle x \rangle_q^n\,\big)  
  &=& n \;\;\: \langle x \rangle_q^{n-1},
 \end{array}
\end{eqnarray}
\begin{eqnarray}
 \label{eq:D-power-nonlinear}
 \begin{array}{llllll}
 \smallskip
 &\widetilde{\text{D}}_{\text{i}} \,\big( \, 
                                     \langle x \rangle_q \ioperator{\owedge} n\,
                                    \big) &=&
 \widetilde{\text{D}}_{\text{i}} \,\big( \;\;     
                                     \langle x^n \rangle_q 
                                   \big) 
  &=& n\,x^{n-1},
\\
 &\widetilde{\text{D}}_{\text{o}} \,\big( \, 
                                   {}_q\langle x \rangle \ooperator{\owedge} n\,
                                    \big) &=&
 \widetilde{\text{D}}_{\text{o}} \,\big(  
                                        {}_q\langle x^n \rangle \,
                                   \big) 
  &=& n\,x^{n-1}.
 \end{array}
\end{eqnarray}

Second and higher deformed linear derivatives follow the usual rule,
$\text{D}_{\delta}^2 f(x) = \text{D}_{\delta} 
                              \big[\text{D}_{\delta}f(x)\big]$ 
and so on, 
but for the deformed nonlinear cases, second order derivatives
(and similarly for higher order derivatives) are defined as
\begin{equation}
 \label{eq:2nd-nonlinear-derivative}
 \widetilde{\text{D}}_{\delta}^2 \, f(x) = h_{\delta}\, 
\big(f(x)\big) \, \frac{\text{d}}{\text{d}x}
 \left[ h_{\delta}\big(f(x)\big) \, \frac{\text{d} f(x)}{\text{d}x} \right].
\end{equation}

The product rule for the deformed linear derivatives 
is identical to the usual one,
$\text{D}_{\delta}\big(f(x)\,g(x)\big) = 
 \text{D}_{\delta}\big(f(x)\big) \, g(x) 
 + f(x) \, \text{D}_{\delta}\big(g(x)\big)$.
The product rule for the deformed nonlinear derivatives is 
\begin{equation}
 \label{eq:product-rule-nonlinear}
 \frac{1}{h_{\delta}\big(f(x) g(x)\big)} \, 
 \widetilde{\text{D}}_{\delta} \bigg( f(x) \, g(x) \bigg)
 =
   \left( 
         \frac{1}{h_{\delta}\big(f(x)\big)} \, 
         \widetilde{\text{D}}_{\delta}\,f(x) 
   \right) \, g(x)
   + f(x) \, 
     \left(
           \frac{1}{h_{\delta}\big(g(x)\big)} \, 
           \widetilde{\text{D}}_{\delta}\,g(x)
     \right).
\end{equation}

The deformed antiderivatives, or indefinite deformed integrals, 
associated to the linear deformed derivatives 
are defined by
\begin{eqnarray}
 \label{eq:q-linear-integral}
 \smallskip
 \int_{(\delta)}^x f(x') \, \text{d}x' &\equiv& \int^x f(x')\,\text{d}_\delta x',
 \\
                                &=&    \int^x f(x')\,h_{\delta}(x')\,\text{d}x'
\end{eqnarray}
(The symbol $(\delta)$ within parenthesis refers to the deformation, 
 and not a limit of integration),
so
\begin{equation}
 \label{eq:q-linear-derivative-of-q-integral}
 \text{D}_\delta \, \int_{(\delta)}^x f(x') \, \text{d}x' = f(x)
\end{equation}
and
\begin{equation}
 \label{eq:q-integral-of-q-linear-derivative}
 \int_{(\delta)}^x \text{D}_\delta f(x') \, \text{d}x' = f(x) + C.
\end{equation}
One possibility for defining the deformed antiderivatives
associated to the nonlinear deformed derivatives,
particularly following the definition used in \cite{Borges-2004}
for the $\delta=$ ole case, is
\begin{eqnarray}
 \label{eq:q-nonlinear-integral}
 \widetilde{\int_{(\delta)}^x} f(x') \, \text{d}x'
   \equiv \int^x \frac{1}{h_{\delta}\big(f(x')\big)} \, f(x') \, \text{d} x',
\end{eqnarray}
A significant weakness with this option is that 
the following important properties are not satisfied:
\begin{equation}
 \label{eq:q-nonlinear-derivative-of-q-integral}
 \widetilde{\text{D}}_\delta \, \widetilde{\int_{(\delta)}^x} \, f(x') \, \text{d}x' 
 \ne f(x)
\end{equation}
and
\begin{equation}
 \label{eq:q-integral-of-q-nonlinear-derivative}
 \widetilde{\int_{(\delta)}^x} \widetilde{\text{D}}_\delta \, f(x') \, \text{d}x' 
 \ne f(x) + C.
\end{equation}

\section{\label{sec:entropy-generator}Entropy generator}

Connection between entropies and derivatives was pointed out by
Abe \cite{Abe-pla-1997}.  
He observed that the Boltzmann-Gibbs entropy can be rewritten as (with $k = 1$)
\begin{equation}
 \label{eq:S_1-generator}
 S_1 = - \left.\frac{\text{d}}{\text{d}\alpha} g(\alpha)\right|_{\alpha=1}
\end{equation}
with
\begin{equation}
 \label{eq:generating-function}
 g(\alpha) = \sum_i^W p^\alpha.
\end{equation}
He realized that $S_q$ entropy can be similarly recast through 
the Jackson's derivative of a function $f(x)$
\cite{Jackson-1909}
\begin{equation}
 \label{eq:jackson}
 \text{D}_q^{(\text{J})} f(x) \equiv \frac{f(qx) - f(x)}{qx - x}
\end{equation}
(the same deformed derivative of quantum calculus \cite{Kac-Cheung-2002};
 Newtonian derivative is recovered as the limiting case $q \to 1$),
so
\begin{equation}
 \label{eq:Sq-jackson}
 S_q = - \left.
               \text{D}_q^{(\text{J})} g(\alpha)
         \right|_{\alpha=1}.
\end{equation}
This property has been interpreted as expressing the association between
Boltzmann-Gibbs entropy ($S_1$) to infinitesimal translations,
and Tsallis entropy to finite dilations
\cite{Tsallis-springer-2009}.
Abe applied this procedure a step further, 
and used a different derivative operator on $g(\alpha)$,
generating a new symmetric entropic functional $S_q^S$ 
with $q \leftrightarrow q^{-1}$ invariance.
Following the same line, a two-parameter derivative operator was used
to define a two-parameter $S_{q,q'}$ entropy, that recovers the previous
$S_q^S$, $S_q$ and $S_1$ with convenient choices of the indices $q$ and $q'$
\cite{Borges-Roditi-1998}.

All the eight deformed derivatives
(\ref{eq:ile-derivatives})--(\ref{eq:oel-derivatives})
applied on (\ref{eq:generating-function}) 
result $S_1$ entropy with a multiplying function of the parameter $q$:
$-\left.\text{D}_{\delta} g(\alpha)\right|_{\alpha=1} 
  = h_{\delta}^\eta(1) \, S_1$,
where $\text{D}_{\delta}$ represents any of the deformed 
(linear or nonlinear) derivatives
(at this point we don't use the tilde for the nonlinear deformed derivatives),
$h_{\delta}(1)$ is a particular value of the corresponding 
Eq.\ (\ref{eq:deforming-functions}), 
$\eta=-1$ for the linear    deformed derivatives, and
$\eta=+1$ for the nonlinear deformed derivatives.
This is consequence of the generalized derivatives
being based on infinitesimal deformed translations, 
and the infinitesimal nature of the translation determines the entropy
(except for a multiplicative constant),
despite of the deformations.

A non-trivial result is obtained by inverting the procedure.
Instead of applying one of the generalized derivatives on the generating function
(\ref{eq:generating-function}),
we apply the ordinary Newtonian derivative on a generalized generating function:
\begin{equation}
 \label{eq:S_q-generator}
 S_q^{\delta} = - \left.
                        \frac{\text{d}}{\text{d}\alpha} g_{\delta}(\alpha;q)
                  \right|_{\alpha=1}.
\end{equation}
The generalized generating functions are obtained through 
the four generalized powers, 
(\ref{eq:ile-power}),
(\ref{eq:ole-power}),
(\ref{eq:iel-power}),
(\ref{eq:oel-power}):
$ g_{\text{ile}}(\alpha;q) 
   = \sum_i^W (p_i {\ileoperator{\owedge} \alpha}) $, \ 
$ g_{\text{ole}}(\alpha;q) 
   = \sum_i^W (p_i {\oleoperator{\owedge} \alpha}) $, \ 
$ g_{\text{iel}}(\alpha;q) 
   = \sum_i^W (p_i {\ieloperator{\owedge} \alpha}) $, \ 
$ g_{\text{oel}}(\alpha;q) 
   = \sum_i^W (p_i {\oeloperator{\owedge} \alpha}) $.
The resulting functionals are
\begin{subequations}
\label{eq:S}
\begin{equation}
 \label{eq:Sile}
 S_q^{\text{ile}} = \sum_i {}_q[-p_i] \ln \big({}_q[p_i]\,\big),
\end{equation}
\begin{equation}
 \label{eq:Sole}
 S_q^{\text{ole}} = -\sum_i [p_i]_q \, \ln\big(\, [p_i]_q \,\big)
                    - (1-q) \sum_i p_i \, [p_i]_q \, \ln\big(\, [p_i]_q \,\big),
\end{equation}
\begin{equation}
 \label{eq:Siel}
 S_q^{\text{iel}} = -\sum_i p_i \ln \big(\,{}_q\{p_i\}\,\big)
                    - (1-q) \sum_i p_i \ln p_i \, \ln \big(\,{}_q\{p_i\}\,\big),
\end{equation}
\begin{equation}
 \label{eq:Soel}
 S_q^{\text{oel}} = 
                   - \displaystyle \sum_i p_i^q \ln \big(\,\{p_i\}_q\,\big).
\end{equation}
\end{subequations}%
          The use of the generalized derivatives essentially results the same,
          $-\left.\text{D}_{\delta} g_{\delta}(\alpha;q)\right|_{\alpha=1} 
           = h_{\delta}^\eta(1) \, S_q^{\delta}$,
          except for a multiplicative constant for the le cases,
          since 
          $h_{\text{iel}}(1)=h_{\text{oel}}(1)=1$.
The certainty distribution originates non zero values for the le functionals:
$ S_q^{\text{ile}}[p_i=1; p_{j}=0, \forall j \ne i] \ne 0 $  
for $q>1$,
and, 
$ S_q^{\text{ole}} [p_i=1; p_{j}=0, \forall j \ne i] \ne 0 $  
for $q<1$,
since ${}_q[1] \ne 1$ and $[1]_q \ne 1$. 
Also, 
the le functionals present negative values:
$S_q^{\text{ile}}$ presents negative values for $q<1$,
$S_q^{\text{ole}}$ presents negative values for $q>1$.
Besides, there are ranges of values of $q$ for which neither
$ S_q^{\text{ile}}$
nor
$ S_q^{\text{ole}}$
present a definite concavity 
(two instances: $q=2.4$, for ile; $q=2.3$, for ole).
These are severe drawbacks and consequently 
(\ref{eq:Sile}) and (\ref{eq:Sole}) 
can not be considered as legitimate entropic forms.

The iel-functional $S_q^{\text{iel}}$ fails on the expansibility 
property for $q<1$ (adding events of zero probability), 
since ${}_{q<1}\{0\}$ is not defined.
For $q>1$, it is expansible, non negative and the certainty distribution 
($p_i=1; p_{j}=0, \forall j \ne i$) implies $S_q^{\text{iel}} = 0$, 
so, (\ref{eq:Siel}) is admissible as an entropic form for $q>1$.

\begin{figure}[!htb]
\centering
 \begin{center}
 \begin{minipage}{0.47\linewidth}
  \includegraphics[width=0.8\linewidth,keepaspectratio,clip=]{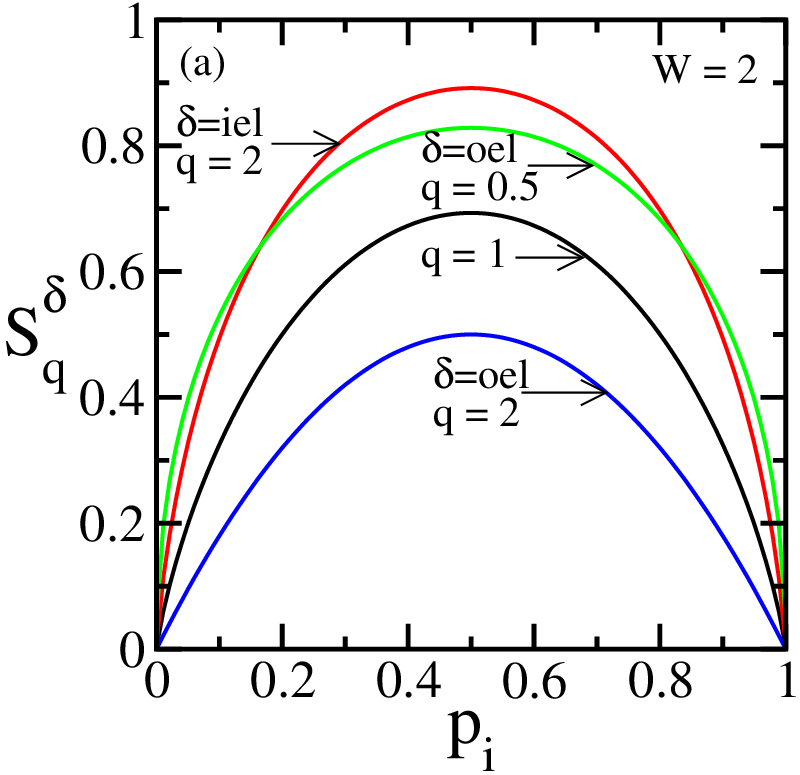}
 \end{minipage}
 \begin{minipage}{0.47\linewidth}
  \includegraphics[width=0.8\linewidth,keepaspectratio,clip=]{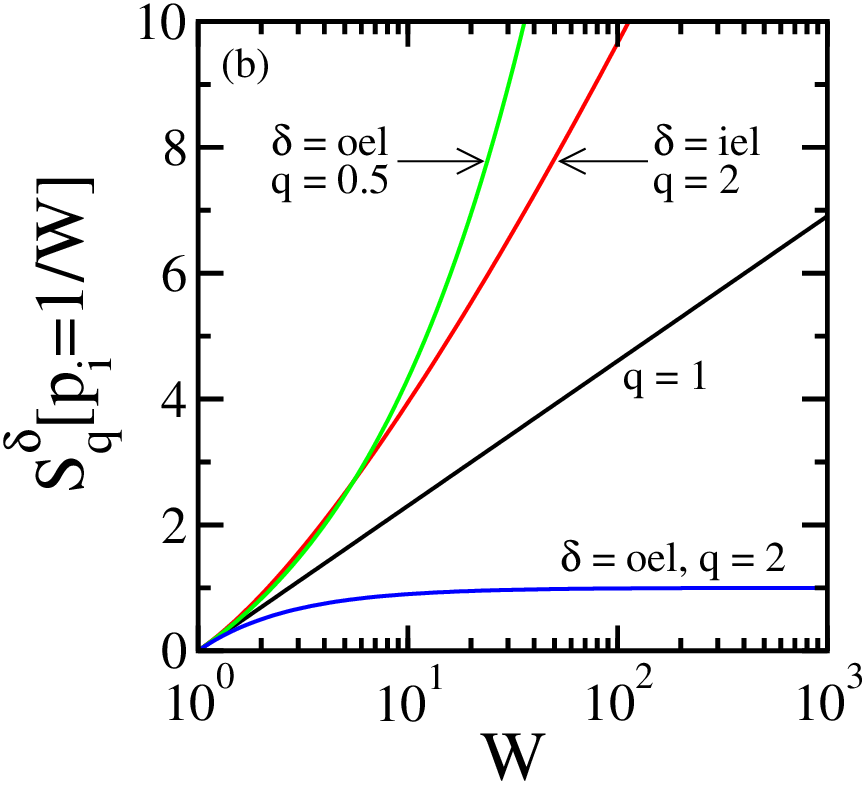}
 \end{minipage}
 \end{center}
 \caption{
          (a) el-Entropies for a two-state system.
          $S_q^{\text{iel}}$ (\protect\ref{eq:Siel})
          for $q=2$ (red);
          $S_q^{\text{oel}}=S_q$ 
          (\protect\ref{eq:Soel}),
          for $q=0.5$ (green), $q=2$ (blue);
          $S_1$ (black).
          $S_q$ entropy is convex for $q<0$,
          see \protect\cite{Tsallis-1988}.
          (b) el-Entropies for equiprobable states as a function of $W$.
          Abscissa in log scale, for which the Boltzmann case 
          is a straight line (black).
          $S_q^{\text{iel}}$ for $q=2$ (red),
          $S_q^{\text{oel}}$ for $q=0.5$ (green), $q=2$ (blue).
          \label{fig:sq}
         }
\end{figure}

The oel-functional (\ref{eq:Soel}) is the nonadditive entropy $S_q$ 
[see Eq.\ (\ref{eq:qlog-as-log})], vastly considered in the literature.
This result permits to amend a previous statement:
$S_q$ entropy, that is associated to finite dilations,
can also be associated to infinitesimal translations,
but in a deformed space expressed by the oel-power.
Figure \ref{fig:sq}a illustrates the concavity for the two
admissible entropic functionals, 
Eq.\ (\ref{eq:Siel}) with $q>1$ and Eq.\ (\ref{eq:Soel}),
for a two-state system.
Figure \ref{fig:sq}b illustrates el-entropies as 
monotonically increasing functions of the number of states $W$
for the equiprobable distribution, $p_i=1/W$, $\forall i$, 
with the abscissa in logarithm scale, 
for which the usual case appears as a straight line.

\section{\label{sec:conclusions}Final remarks}

A forerunner of the transformations given by Eq.\ (\ref{eq:le-numbers})
is the relation between R\'enyi entropy,
$S^{\text{R}}_q = (1-q)^{-1} \ln\big(\sum_i^W p_i^q\big)$,
and Tsallis entropy (\ref{eq:Sq})
[see Eq.\ (8) of Ref.\ \cite{Tsallis-1988}],
$S_q^{\text{R}} = [S_q]_q$,
and, equivalently,
$S_q = {}_q[S_q^{\text{R}}]$.
Another instance of the transformation
represented by the ile-number (\ref{eq:ile-number})
appeared in Eq.\ (22) of Ref.\ \cite{Borges-1998}
and allowed the generalization of trigonometric functions.
The ole-number ${}_q[x]$ appeared as Eq.\ (5) 
of Ref.\ \cite{Czachor-Naudts-2002}, 
as the scaling factor of the generalized Kolmogorov-Nagumo average 
for expressing the R\'enyi entropy.
A former example of connecting deformed numbers with 
deformed differential operators have appeared in
Ref.\ \cite{daCosta-Borges-2014,daCosta-Borges-2018},
with
the transformation (\ref{eq:ile-number})
and the deformed differential (\ref{eq:linear-ole-derivative}),
establishing an equivalence between a position-dependent mass system 
in a usual space
and a constant mass within a deformed space.
These works have been recently extended to the deformed version 
of the Fokker-Planck equation for inhomogeneous medium 
with position-dependent mass \cite{daCosta-Gomez-Borges-2020}.

Expressions with operations belonging to one class of $q$-algebra
may result in operations belonging to a different class.
Some instances: the following are generalizations of the logarithm of a product 
as a sum of logarithms \big($\ln(xy)=\ln x + \ln y$\big):
\begin{subequations}
 \label{eq:log-of-product}
 \begin{equation}
  \ln_q\big(\,xy\,\big) 
  = \ln_q x \oleoperator{\oplus} \ln_q y,
 \end{equation}
 \begin{equation}
  \ln_q\big( x \oeloperator{\otimes} y \, \big) 
  = \ln_q x + \ln_q y,
 \end{equation}
 \begin{equation}
  \ln  \big( x \ieloperator{\otimes} y \, \big) 
  = \ln x \oleoperator{\oplus} \ln y,
 \end{equation}
 \begin{equation}
  \ln  \big( x \oeloperator{\otimes} y \, \big) 
  = \ln x \ileoperator{\oplus} \ln y.
 \end{equation}
\end{subequations}%
Generalizations of the logarithm of a power, $\ln x^y = y \ln x$, are
\begin{subequations}
 \label{eq:log-of-power}
 \begin{equation}
  \label{eq:qlog-of-oelpower}
  \ln_q \big( x \oeloperator{\owedge} y \, \big) = y \ln_q x. 
 \end{equation}
 \begin{equation}
  \ln \big( x \ieloperator{\owedge} y \, \big) = y \oleoperator{\odot} \ln x, 
 \end{equation}
 \begin{equation}
  \ln \big( x \oeloperator{\owedge} y \, \big) = y \ileoperator{\odot} \ln x. 
 \end{equation}
\end{subequations}%
The counterpart of these expressions are generalizations of
the exponential of a sum as a product of exponentials,
$\text{e}^{x+y} = \text{e}^x \, \text{e}^y$:
\begin{subequations}
 \label{eq:exp-of-sum}
 \begin{equation}
  \text{e}_q^{\,x \oleoperator{\oplus} y}  
  = \text{e}_q^{\,x} \: \text{e}_q^{\,y},
 \end{equation}
 \begin{equation}
  \text{e}_q^{\, x \,+\, y} 
  = \text{e}_q^{\,x} \oeloperator{\otimes} \text{e}_q^{\,y},
 \end{equation}
 \begin{equation}
  \text{e}^{\: x \oleoperator{\oplus} y}  
  = \text{e}^{\,x} \ieloperator{\otimes} \text{e}^{\,y},
 \end{equation}
 \begin{equation}
  \text{e}^{\: x \ileoperator{\oplus} y } 
  = \text{e}^{\,x} \oeloperator{\otimes} \text{e}^{\,y},
 \end{equation}
\end{subequations}%
and the power of an exponential as the exponential of a product
\big($(\text{e}^x)^{\,y}=\text{e}^{yx}$\big):
\begin{subequations}
 \label{eq:power-of-exp}
 \begin{equation}
 \text{e}_q^{\,x} \oeloperator{\owedge} y = \text{e}_q^{\: y x },
 \end{equation}
 \begin{equation}
 \text{e}^{\,x} \ieloperator{\owedge} y =\text{e}^{\: y \oleoperator{\odot} x },
 \end{equation}
 \begin{equation}
 \text{e}^{\,x} \oeloperator{\owedge} y =\text{e}^{\: y \ileoperator{\odot} x }.
 \end{equation}
\end{subequations}%
Relations~(\ref{eq:log-of-product}) are also valid 
for the logarithm, or for the $q$-logarithm, of a ratio,
simply replacing ordinary or general products by ordinary or general ratios,
and ordinary or general sums by ordinary or general differences.
Similarly, relations~(\ref{eq:exp-of-sum}) are also valid 
for the exponential, or for the $q$-exponential, of a difference,
by replacing the operators accordingly.

Equation (\ref{eq:qlog-of-oelpower}) is the one referred to in the Introduction,
that makes $S_q$ extensive:
consider a composed system for which its subsystems have $W_i > 1$ available
states. If they are independent, the number of available states of the
composed system is $W=\prod_i^N W_i$, and, 
besides, if they are identical, $W=W_1^N$.
Correlations between the subsystems lead to a smaller number of available states
for the composed system, and particular strong correlations represented by
$W = W_1 \oeloperator{\owedge}N$, with $q<1$,
makes $S_q = k \ln_q W = N k \ln_q W_1$.
This is a non trivial case of extensivity.

Different possibilities for generating rules of arithmetic operations,
instead of (\ref{eq:generating-rules}), are
${}_q[x]   \oleoperator{\Circle} {}_q[y]   = {}_q[ x \circ y ]$,
\,
$    [x]_q \ileoperator{\Circle}     [y]_q =     [ x \circ y ]_q$;
Ref.\ \cite{Lobao-bjp-2009} uses these patterns.

Weberszpil, Lazo and Helay\"el-Neto \cite{Weberszpil-2015}
have shown that the linear ole-derivative 
(\ref{eq:linear-ole-derivative})
is the first order expansion of the Hausdorff derivative. 
Whether the other generalized derivatives are also connected 
to fractal derivatives and fractal metrics remains to be investigated.

Two of the functionals obtained with the recipe of applying the ordinary
derivatives to a generalized version of the generating function,
(\ref{eq:S_q-generator}),
result admissible entropic forms corresponding to the el-class:
$S_q^{\text{iel}}$ (\ref{eq:Siel}), 
and 
$S_q^{\text{oel}}$ (\ref{eq:Soel}).
The other functionals (\ref{eq:Sile}) and (\ref{eq:Sole})
are not admissible to be considered as entropies,
but this does not mean that the le-algebras or le-calculus
they are based on are not feasible for other applications.

Extension to the complex domain of the deformed numbers 
still remains to be explored.
Two-parameter generalization are not addressed here, 
we just advance a few lines.
Two-parameter generalizations of numbers in accordance with the present 
developments are given by
\begin{subequations}
  \label{eq:q,q'-numbers}
\begin{eqnarray}
  \label{eq:qq'-le-number}
   {}_q[x]_{q'} &\equiv& [x]_{q,q'} = \ln_q \exp_{q'}(x),
  \\
  \label{eq:qq'-el-number}
   {}_q\{x\}_{q'} &\equiv& \{x\}_{q,q'}  \; = \; \exp_q(\ln_{q'} x).
\end{eqnarray}
\end{subequations}%
The use of the relatively uncommon subscripted prefix to represent
the two parameter deformed number may be avoided, 
since there is no ambiguity with the symbol
$\langle x\rangle_{q,q'}$.
Two-parameter arithmetic operators follow straightforwardly:
\begin{eqnarray}
 \label{eq:qq'-generating-rule}
  x \; \Circle_{\langle q,q'\rangle} \; y = 
                                      \left\langle \: 
                                       \langle x \rangle_{q',q} \; 
                                      \circ \; 
                                       \langle y \rangle_{q',q} \: 
                                      \right\rangle_{q,q'},
\end{eqnarray}
for which, of course, all the previous developments are particular cases.
The two-parameter algebra of Ref.\ \cite{Cardoso-Borges-Lobao-Pinho-2008}
is obtained through a different generating rule than 
(\ref{eq:qq'-generating-rule});
it derives from the two-parameter generalized
logarithm and exponential functions \cite{Schwammle-Tsallis-2007}
[Eq.\ (16) and (17) of \cite{Cardoso-Borges-Lobao-Pinho-2008}].

It also comes naturally the two-parameter derivative $\text{D}_{q,q'}f(x)$, 
with deformation on both the independent and dependent variables.
A broader generalization of the derivatives can be defined 
by using not only deformations on the variation of the independent 
and dependent variable, but also on the ratio among them,
with three parameters, 
in a rather intricate way, say:
$q$   for the deformed differential of the independent variable, 
$q'$  for the deformed differential of the   dependent variable, and
$q''$ for the deformed ratio between them. 
A particular case with $q=q'=q''$ 
was done in Ref.\ \cite{Kalogeropoulos-2005},
and, more recently in Ref.\ \cite{Czachor-2020}.

Finally, all the present scenario stands on the pair of 
$q$-logarithm/$q$-exponential functions, inverse of each other.
The whole picture may be differently deformed by using different
continuous, monotonous, invertible pair of functions,
in agreement with Ref.\ (\ref{eq:generalized-numbers}).

\section*{Acknowledgments}
%
This work was partially supported by National
Institute of Science and Technology for Complex Systems (INCT-SC).
E.\ P.\ B.\ thanks C.\ Tsallis and Si Hyung Joo,
and both authors thank I.\ S.\ Gomez for stimulating discussions.

\appendix
\section{A note on notations --- Explicit expressions}

The peculiar notation adopted in the present work
is conceived for compactness, 
once the explicit forms of some equations may be large or cumbersome.
The notation for the generalized numbers has been inspired in the
{\it $Q$-analog of $n$} \cite{Kac-Cheung-2002},
a generalized number represented within square brackets 
(\ref{eq:quantum-number}).
The four classes of generalized numbers are grouped into two categories,
one, the `le' category, 
uses the generalized exponential (or its ordinary version) as argument
of the ordinary logarithm (or its generalized version), 
and the other, the `el' category, the other way around.
We have used square brackets for the former, and curly brackets for the later.
Some ambiguity is unfortunately unavoidable, as square and curly brackets
are also used with their usual meanings, 
and the reader must resolve it by the context.
We refer to them as `le' or `el' concerning the order in which the
logarithm/exponential functions appear.  
Despite of the unusualness, or even possibly strangeness, of the notation, 
we consider it may help identify the classes more promptly than 
something like `type 1', `type 2' etc.
Differently from the generalized numbers,
we use the subscripts enclosed by their corresponding brackets,
when dealing with generalized arithmetic operators, 
so the reader can easily identify the object being generalized,
if it is a number or an operator.
We have chosen prefix and postfix subscripts, to avoid using superscripts.
These pair of subscripts may play a simplifying role if used appropriately,
as illustrated by Eq.\ (\ref{eq:inverse}).
In the following we present explicit forms of some expressions,
for the benefit of the interested reader.
The notation $[\cdot]_+ \equiv \text{max}\{0,\cdot\}$ is still used here.

\bigskip

\noindent
{\bf ile-number} [Eq.\ (\ref{eq:ile-number})]
\begin{equation}
 \label{eq:ile-number-explicit}
 [x]_q = \frac{1}{1-q} \ln\bigg(1+(1-q)x\bigg)_+.
\end{equation}

\medskip

\noindent
{\bf ole-number} [Eq.\ (\ref{eq:ole-number})]
\begin{equation}
 \label{eq:ole-number-explicit}
 {}_q[x] = \frac{\text{e}^{(1-q)x}-1}{1-q}.
\end{equation}

\medskip

\noindent
{\bf iel-number} [Eq.\ (\ref{eq:iel-number})]
\begin{equation}
 \label{eq:iel-number-explicit}
 \{x\}_q = \text{sign}(x) \: \exp\left(\frac{|x|^{1-q}-1}{1-q}\right).
\end{equation}

\medskip

\noindent
{\bf oel-number} [Eq.\ (\ref{eq:oel-number})]
\begin{equation}
 \label{eq:oel-number-explicit}
 {}_q\{x\} = \text{sign}(x) \: \bigg(1+(1-q)\ln |x|\bigg)_+^{1/(1-q)}.
\end{equation}

\medskip

\noindent
{\bf ile-addition, ile-subtraction}
     [Eq.\ (\ref{eq:ile-sum}), (\ref{eq:ile-difference})]
\begin{eqnarray}
 \label{eq:ile-pm-explicit}
 x \ileoperator{\opm} y &=&
 \frac{1}{1-q} \ln 
                \left[
                  1+(1-q)\left(
                          \frac{\text{e}^{(1-q)x}-1}{1-q}
                          \pm
                          \frac{\text{e}^{(1-q)y}-1}{1-q}
                         \right)
                \right]_+,
  \\
  \bigskip
                          &=&
 \frac{1}{1-q} \ln 
                \left[
                 \text{e}^{(1-q)x} \pm \text{e}^{(1-q)y} \mp 1
                \right]_+.
\end{eqnarray}

\medskip

\noindent
{\bf ile-multiplication} [Eq.\ (\ref{eq:ile-product})]
\begin{equation}
 \label{eq:ile-product-explicit}
 x \ileoperator{\otimes} y =
 \frac{1}{1-q} \ln 
                \left[
                  1+\frac{(\text{e}^{(1-q)x}-1) \: (\text{e}^{(1-q)y}-1)}{1-q}
                \right]_+.
\end{equation}

\noindent
{\bf ile-division} [Eq.\ (\ref{eq:ile-ratio})]
\begin{equation}
 \label{eq:ile-ratio-explicit}
 x \ileoperator{\oslash} y =
 \frac{1}{1-q} \ln 
                \left[
                  1+(1-q) \, \frac{\text{e}^{(1-q)x}-1}{\text{e}^{(1-q)y}-1}
                \right]_+.
\end{equation}

\medskip

\noindent
{\bf ile-power} [Eq.\ (\ref{eq:ile-power})]
\begin{equation}
 \label{eq:ile-power-explicit}
 x \ileoperator{\owedge} y =
 \frac{1}{1-q} \ln 
                \left[
                  1+(1-q) \,\left( \frac{\text{e}^{(1-q)x}-1}{1-q} \right)^y
                \right]_+,
                \quad (x>0).
\end{equation}

\noindent
{\bf ole-addition} [see Eq.\ (\ref{eq:ole-sum})]

\medskip

\noindent
{\bf ole-subtraction} [see Eq.\ (\ref{eq:ole-difference})]

\medskip

\noindent
{\bf ole-multiplication} [Eq.\ (\ref{eq:ole-product})]
\begin{eqnarray}
 \label{eq:ole-product-explicit}
  x \oleoperator{\otimes} y &=&
     \frac{
           \exp\left[ 
                 \frac{\ln\left[1+(1-q)x\right]_+  \ln\left[1+(1-q)y\right]_+}
                      {1-q}
               \right]
               - 1
          }
          {1-q},
  \\
  \bigskip
                            &=&
     \frac{
           \left[1+(1-q)x \right]_+^{\frac{1}{1-q}\ln \left[ 1+(1-q)y \right]_+}
            - 1
          }
          {1-q},
  \\
  \bigskip
                            &=&
     \frac{
           \left[1+(1-q)y \right]_+^{\frac{1}{1-q}\ln \left[ 1+(1-q)x \right]_+}
            - 1
          }
          {1-q}.
\end{eqnarray}

\medskip

\noindent
{\bf ole-division} [Eq.\ (\ref{eq:ole-ratio})]
\begin{equation}
 \label{eq:ole-ratio-explicit}
  x \oleoperator{\oslash} y =
     \frac{
           \exp\left[ (1-q)
                 \frac{\ln\left[1+(1-q)x\right]_+}{\ln\left[1+(1-q)y\right]_+}
               \right]
               - 1
          }
          {1-q}.
\end{equation}

\medskip

\noindent
{\bf ole-power} [Eq.\ (\ref{eq:ole-power})]
\begin{equation}
 \label{eq:ole-power-explicit}
  x \oleoperator{\owedge} y =
     \frac{
           \exp\left[ (1-q)^{1-y}
                      \ln^y\left[1+(1-q)x\right]_+
               \right]
               - 1
          }
          {1-q},
          \quad (x>0).
\end{equation}

\medskip

\noindent
{\bf iel-addition, iel-subtraction}
     [Eq.\ (\ref{eq:iel-sum}), (\ref{eq:iel-difference})]
\begin{eqnarray}
 \label{eq:iel-pm-explicit}
 x \ieloperator{\opm} y 
   &=&
   \text{sign}(x \pm y) \,
 \nonumber \\
   &\times&
   \exp\left(
        \frac{
              \left|
                    \;
                    \text{sign}(x)\,\bigg[
                                            1+(1-q)\ln |x|
                                      \bigg]_+^{\frac{1}{1-q}}
                    \pm \;
                    \text{sign}(y)\,\bigg[
                                            1+(1-q)\ln |y|
                                       \bigg]_+^{\frac{1}{1-q}}
                    \;
              \right|^{1-q}
              - 1 }
             {1-q}
       \right).
\end{eqnarray}

\medskip

\noindent
{\bf iel-multiplication}
     [Eq.\ (\ref{eq:iel-product})]
\begin{equation}
 \label{eq:iel-product-explicit}
 x \ieloperator{\otimes} y 
   =
   \text{sign}(xy) \,
   \exp\left(
         \frac{
               \left|
                     \;
                     \bigg[1+(1-q) \ln |x|\bigg]_+
                     \;
                     \bigg[1+(1-q) \ln |y|\bigg]_+
                     \;
               \right|
               - 1 }
              {1-q}
       \right).
\end{equation}

\medskip

\noindent
{\bf iel-division}
     [Eq.\ (\ref{eq:iel-ratio})]
\begin{equation}
 \label{eq:iel-ratio-explicit}
 x \ieloperator{\oslash} y 
   =
   \text{sign}(x/y) \,
   \exp\left[
              (1-q)^{-1} \;
          \left( \;
               {\left|
                     \;
                     \frac{\bigg[1+(1-q) \ln |x|\bigg]_+}
                          {\bigg[1+(1-q) \ln |y|\bigg]_+}
                     \;
               \right|
               - 1 }
          \right) \;
       \right].
\end{equation}

\medskip

\noindent
{\bf iel-power} [Eq.\ (\ref{eq:iel-power})]
\begin{equation}
 \label{eq:iel-power-explicit}
 x \ieloperator{\owedge} y 
   =
   \text{sign}(x) \,
   \exp\left(
         \frac{
               \left|
                     \;
                     \bigg[1+(1-q) \ln |x|\bigg]_+^y
                     \;
               \right|
               - 1 }
              {1-q}
       \right),
       \quad (x>0).
\end{equation}

\medskip

\noindent
{\bf oel-addition, oel-subtraction}
     [Eq.\ (\ref{eq:oel-sum}), (\ref{eq:oel-difference})]
\begin{eqnarray}
 \label{eq:oel-pm-explicit}
 x \oeloperator{\opm} y 
   &=&
   \text{sign}(x \pm y) \,
 \nonumber \\
   &\times&
   \left[
    1 + (1-q) \ln \left|\:
                        \text{sign}(x)\,
                        \exp\left(
                                  \frac{|x|^{1-q}-1}{1-q}
                             \right)
                        \; \pm \;
                        \text{sign}(y)\,
                        \exp\left(
                                  \frac{|y|^{1-q}-1}{1-q}
                             \right) \:
                  \right| \;
   \right]_+^{\frac{1}{1-q}}.
\end{eqnarray}

\medskip

\noindent
{\bf oel-multiplication} [see Eq.\ (\ref{eq:oel-product})]

\medskip

\noindent
{\bf oel-division} [see Eq.\ (\ref{eq:oel-ratio})]

\medskip

\noindent
{\bf oel-power} [Eq.\ (\ref{eq:oel-power})]
\begin{eqnarray}
 \label{eq:oel-power-explicit}
  x {\oeloperator{\owedge} y} = \big(\text{sign}(x)\big)^y \,
                                \left[
                                      y\,|x|^{1-q} - (y-1)
                                \right]_+^{\frac{1}{1-q}},
                                \quad (x>0).
\end{eqnarray}

\medskip

\noindent
$S_q^{\text{ile}}$ {\bf functional} [Eq.\ (\ref{eq:Sile})]
\begin{eqnarray}
 \label{eq:Sile-explicit}
 S_q^{\text{ile}} 
 = \sum_i^W \frac{\text{e}^{-(1-q)p_i}-1}{1-q} 
   \ln\left[
            \frac{\text{e}^{(1-q)p_i}-1}{1-q} 
      \right].
\end{eqnarray}

\medskip

\noindent
$S_q^{\text{ole}}$ {\bf functional} [Eq.\ (\ref{eq:Sole})]
\begin{eqnarray}
 \label{eq:Sole-explicit}
 S_q^{\text{ole}} 
 = - \sum_i^W 
               \frac{\ln\big[1+(1-q)p_i\big]}{1-q}
                            \ln \left[
                                      \frac{\ln\big[1+(1-q)p_i\big]}{1-q}
                                \right]
\nonumber \\
   - \sum_i^W  
               p_i \ln\big[1+(1-q)p_i\big]
                   \ln \left[
                             \frac{\ln\big[1+(1-q)p_i\big]}{1-q}
                       \right]
\end{eqnarray}

\medskip

\noindent
$S_q^{\text{iel}}$ {\bf functional} [Eq.\ (\ref{eq:Siel})]
\begin{eqnarray}
 \label{eq:Siel-explicit}
 S_q^{\text{iel}} 
 = - \sum_i^W \frac{p_i}{1-q} \ln\big[1+(1-q) \ln p_i\big]
   - \sum_i^W p_i \ln p_i \ln\big[1+(1-q)] \ln p_i\big]
\end{eqnarray}

\medskip

\noindent
$S_q^{\text{oel}}$ {\bf functional} 
 [Eq.\ (\ref{eq:Soel}), see also Eq.\ (\ref{eq:Sq})]

\section*{References}


\end{document}